\documentclass[11pt, reqno, psamsfonts]{amsart}
\pdfoutput=1

\usepackage{amssymb}
\usepackage{amsthm}
\usepackage{amsmath}
\usepackage{latexsym}
\usepackage[T1]{fontenc}
\usepackage[utf8]{inputenc}
\usepackage[russian, french, english]{babel}

\usepackage{graphicx}
\usepackage{wrapfig}
\usepackage[justification=centering, labelfont=bf]{caption}
\usepackage{subcaption}
\usepackage{mathtools}
\usepackage[hidelinks]{hyperref}
\usepackage{tikz}
\usepackage{amsbsy}
\usepackage[inline]{enumitem}
\usepackage{mathrsfs}
\usetikzlibrary{shapes,snakes}
\usetikzlibrary{arrows.meta}
\usetikzlibrary{decorations.pathreplacing}
\usepackage{float}
\usepackage{multicol}
\usepackage{stmaryrd}
\usepackage{eucal}

\usepackage{caladea}
\usepackage[T1]{fontenc}
\usepackage[euler-digits,small]{eulervm}
\usepackage[bb=pazo,scr=esstix]{mathalfa}

\usepackage{upgreek}
\usepackage{titlesec}
\usepackage[spacing=true,kerning=true,babel=true,tracking=true]{microtype}
\usepackage[shortcuts]{extdash}
\usepackage[foot]{amsaddr} 
\usepackage{mdframed}
\usepackage{algpseudocode}
\usepackage{algorithm}
\definecolor{shadecolor}{gray}{0.95}
\usepackage[left=1in,right=1in,top=1in,bottom=1in,bindingoffset=0cm]{geometry}
\usepackage{xspace}
\usepackage{bm}

\usepackage[backend=biber, style=alphabetic, sorting=nyt, maxnames=100,backref=true]{biblatex}
\title{\lsstyle A Fast Distributed Algorithm for $(\Delta + 1)$-Edge-Coloring}
\date{}
\author{\lsstyle Anton~Bernshteyn}
\address{\textls{\normalfont{}Department of Mathematical Sciences, Carnegie Mellon University, Pittsburgh, PA, USA} and \textls{\normalfont{}School of Mathematics, Georgia Institute of Technology, Atlanta, GA, USA}}
\email{bahtoh@gatech.edu}
\thanks{This research was partially supported by the NSF grant DMS-2045412.}

\newtheoremstyle{bfnote}%
{}{}%
{\slshape}{}%
{\bfseries}{\bfseries.}%
{ }%
{\thmname{#1}\thmnumber{ #2}\thmnote{ \ep{\normalfont{}#3}}}

\newtheoremstyle{defbfnote}%
{}{}%
{}{}%
{\bfseries}{.}%
{ }%
{\thmname{#1}\thmnumber{ #2}\thmnote{ (#3)}}

\newtheoremstyle{claim}%
{}{}%
{\slshape}{}%
{\itshape}{.}%
{ }%
{\thmname{#1}\thmnumber{ #2}\thmnote{ \ep{\normalfont{}#3}}}

\theoremstyle{bfnote}

\newtheorem{theo}[equation]{Theorem}

\newtheorem{lemma}[equation]{Lemma}
\newtheorem{corl}[equation]{Corollary}

\newtheorem{fact}[equation]{Fact}

\newtheorem*{claim*}{Claim}

\newtheorem*{corl*}{Corollary}

\newtheorem{case}{Case}
\newtheorem{subcase}{Subcase}[case]
\newtheorem{subsubcase}{Subsubcase}[subcase]

\theoremstyle{claim}

\newcounter{ForClaims}[section]

\newcommand*{\myproofname}{Proof}
\newenvironment{claimproof}[1][\myproofname]{\begin{proof}[#1]}{\end{proof}}

\theoremstyle{definition}

\newtheorem*{defn*}{Definition}

\newtheorem{ques}[equation]{Question}

\newtheorem*{exmp*}{Example}

\newtheorem*{assum*}{Assumptions}

\theoremstyle{remark}
\newtheorem*{ques*}{Question}
\newtheorem*{remk*}{Remark}

\newcommand{\0}{\varnothing}
\newcommand{\set}[1]{\{#1\}}

\newcommand{\dom}{\mathrm{dom}}

\newcommand{\N}{{\mathbb{N}}}
\newcommand{\Z}{\mathbb{Z}}
\newcommand{\F}{\mathbb{F}}
\newcommand{\R}{\mathbb{R}}
\renewcommand{\C}{\mathbb{C}}
\newcommand{\Q}{\mathbb{Q}}

\renewcommand{\epsilon}{\varepsilon}
\renewcommand{\phi}{\varphi}
\renewcommand{\theta}{\vartheta}
\renewcommand{\leq}{\leqslant}
\renewcommand{\geq}{\geqslant}

\renewcommand{\G}{\Gamma}

\newcommand{\defeq}{\coloneqq}
\newcommand{\rest}[2]{{{#1}\vert{#2}}}
\newcommand{\emphd}[1]{{\fontseries{b}\selectfont{#1}}}

\renewcommand{\P}{\mathbb{P}}

\newcommand{\LOCAL}{$\mathsf{LOCAL}$\xspace}
\newcommand{\blank}{\mathtt{blank}}
\newcommand{\Shift}{\mathtt{Shift}}
\newcommand{\length}{\mathtt{length}}
\newcommand{\start}{\mathtt{Start}}
\newcommand{\eend}{\mathtt{End}}
\newcommand{\vstart}{\mathtt{vStart}}
\newcommand{\vend}{\mathtt{vEnd}}
\newcommand{\pivot}{\mathtt{Pivot}}

\newcommand{\bemph}[1]{{\normalfont#1}} 
\newcommand{\ep}[1]{\bemph{(}#1\bemph{)}} 

\newenvironment{scproof}[1][Proof]{\begin{proof}[\textsc{#1}]}{\end{proof}}

\numberwithin{equation}{section}

\renewcommand{\thesubsection}{\arabic{section}.\Alph{subsection}}

\usepackage{etoolbox}

\titleformat{\section}[block]{\bfseries\lsstyle\filcenter}{\thesection.}{1ex}{}
\titleformat{\subsection}[block]{\itshape\filcenter}{\bfseries\upshape\thesubsection.}{1ex}{}
\titleformat{\subsubsection}[runin]{\itshape}{\bfseries\upshape\thesubsubsection.}{1ex}{}[.---]

\titlespacing*{\section}{0pt}{*3}{*1}
\titlespacing*{\subsection}{0pt}{*2.7}{*1}
\titlespacing*{\subsubsection}{0pt}{*1.5}{*0}

\makeatletter
\newcommand{\neutralize}[1]{\expandafter\let\csname c@#1\endcsname\count@}
\makeatother

\newenvironment{theocopy}[1]
{\renewcommand{\theequation}{\ref{#1}}%
	\neutralize{equation}\phantomsection
	\begin{theo}}
	{\end{theo}}

\renewcommand{\mkbibnamefirst}[1]{\textsc{#1}}
\renewcommand{\mkbibnamelast}[1]{\textsc{#1}}
\renewcommand{\mkbibnameprefix}[1]{\textsc{#1}}
\renewcommand{\mkbibnameaffix}[1]{\textsc{#1}}
\renewcommand{\finentrypunct}{}

\renewbibmacro{in:}{}

\renewbibmacro*{volume+number+eid}{%
	\printfield{volume}%
	\setunit*{\addnbspace}
	\printfield{number}%
	\setunit{\addcomma\space}%
	\printfield{eid}}

\DeclareFieldFormat[article]{volume}{\textbf{#1}\space}
\DeclareFieldFormat[article]{number}{\mkbibparens{#1}}

\DeclareFieldFormat{journaltitle}{#1,}
\DeclareFieldFormat[thesis]{title}{\mkbibemph{#1}\addperiod}
\DeclareFieldFormat[article, unpublished, thesis]{title}{\mkbibemph{#1},}
\DeclareFieldFormat[book]{title}{\mkbibemph{#1}\addperiod}
\DeclareFieldFormat[unpublished]{howpublished}{#1, }

\DeclareFieldFormat{pages}{#1}

\DeclareFieldFormat[article]{series}{Ser.~#1\addcomma}

\renewcommand*{\bibfont}{\small}

\setlist{topsep=3pt,itemsep=3pt}

\begin{document}
	\pagestyle{plain}
	
	\vspace*{-15pt}
	
	\maketitle
	
	\begin{abstract}
		We present a deterministic distributed algorithm in the \LOCAL model that finds a proper $(\Delta + 1)$-edge-coloring of an $n$-vertex graph of maximum degree $\Delta$ in $\mathrm{poly}(\Delta, \log n)$ rounds. This is the first nontrivial distributed edge-coloring algorithm that uses only $\Delta+1$ colors \ep{matching the bound given by Vizing's theorem}.	Our approach is inspired by the recent proof of the measurable version of Vizing's theorem due to Greb\'\i{}k and Pikhurko.
	\end{abstract}
	
	\section{Introduction}
	
	\subsection{Distributed algorithms and edge-coloring}
	
	The \LOCAL model of distributed computation was introduced by Linial in the seminal paper \cite{Linial}. In this model an $n$-vertex \ep{simple undirected} graph $G$ abstracts a communication network where each vertex plays the role of a processor and edges represent communication links. Every vertex is given a unique $\Theta(\log n)$-bit identifier \ep{used for symmetry\-/breaking}. Initially, each vertex knows its own identifier, as well as $n$ \ep{the number of vertices} and perhaps some other global parameters, such as the maximum degree $\Delta$ of $G$. The computation proceeds in \emph{rounds}. During each round, the vertices first perform arbitrary local computations and then synchronously broadcast messages to all their neighbors. At the end, each vertex should output its part of the global solution \ep{for instance, its own color or, in the context of edge-coloring, the colors of the edges incident to it}. The only measure of efficiency for such an algorithm is the worst-case number of communication rounds. The reader is referred to the book \cite{BE} by Barenboim and Elkin for an introduction to this subject.
	
	It is clear that in a \LOCAL algorithm that terminates in $T$ rounds, each vertex only has access to information in its radius-$T$ neighborhood. Furthermore, every $T$-round \LOCAL algorithm can be transformed into one in which every vertex first collects all the information about its radius-$T$ neighborhood and then makes a decision, based on this information alone, about its part of the output \ep{see \cite[\S{}4.1.2]{BE}}. In particular, any function of $G$ can trivially be computed by a \LOCAL algorithm in $O(\mathrm{diameter}(G))$ rounds. While in many cases much more efficient distributed algorithms are known, there are still some classical problems for which the trivial bound $O(\mathrm{diameter}(G))$ is the state of the art. In this paper we study one such problem: edge-coloring using $\Delta+1$ colors.
	
	Recall that a \emphd{{proper} $k$-edge-coloring} of a graph $G$ is a map $\phi \colon E(G) \to [k]$ such that $\phi(e) \neq \phi(h)$ for every pair of distinct edges $e$, $h \in E(G)$ that share an endpoint. \ep{Here and in what follows we use the standard notation $[k] \defeq \set{1, \ldots, k}$ for $k \in \N^+$.} A celebrated theorem of Vizing \cite{Viz64} \ep{see \cite[Appendix A.1]{StiebitzBook} for an English translation of Vizing's paper} asserts that every graph $G$ of maximum degree $\Delta$ has a proper $(\Delta+1)$-edge-coloring. However, all heretofore known proofs of Vizing's theorem have been inherently ``non-local'' \ep{we discuss this in more detail in \S\ref{subsec:intro_comb}}. In particular, prior to this work, no nontrivial distributed $(\Delta+1)$-edge-coloring algorithm has been known. Here we provide such an algorithm; specifically, we show that a proper $(\Delta+1)$-edge-coloring of an $n$-vertex graph $G$ can be found by a \LOCAL algorithm whose number of rounds is polynomial in $\log n$ and $\Delta$.
	
	\begin{theo}\label{theo:main}
		There exists a deterministic distributed algorithm in the \LOCAL model that computes a proper $(\Delta + 1)$-edge-coloring of an $n$-vertex graph of maximum degree $\Delta$ in $\mathrm{poly}(\Delta, \log n)$ rounds.
	\end{theo}
	
	While Theorem~\ref{theo:main} provides the first efficient distributed algorithm for $(\Delta + 1)$-edge-coloring, significant progress has been made previously concerning edge-coloring with more colors. In particular, considerable work has been done on the problem of $(2\Delta - 1)$-edge-coloring, since $2\Delta - 1$ is the number of colors required for a greedy \ep{sequential} edge-coloring algorithm. This is a special case of the more general $(\Delta + 1)$-\emph{vertex}-coloring problem, which has also attracted a lot of attention. Among the earliest results in distributed computing is a randomized $O(\log n)$-round $(\Delta + 1)$-{vertex}-coloring algorithm due to Alon, Babai, and Itai \cite{ABI} and independently Luby \cite{Luby}. On the other hand, Goldbreg, Plotkin, and Shannon \cite{GPSh} designed a deterministic $(\Delta + 1)$-vertex-coloring algorithm that runs in $O(\Delta^2) + \log^\ast n$ rounds \ep{here $\log^\ast n$ is the iterated logarithm of $n$, i.e., the number of times the logarithm function must be applied to $n$ before the result becomes at most $1$}. The Goldbreg--Plotkin--Shannon algorithm is extremely fast when $\Delta$ is small compared to $n$ \ep{e.g., when $\Delta$ is constant}, but for large $\Delta$ it becomes slower than the Alon--Babai--Itai/Luby algorithm. Developing a \emph{deterministic} $(\Delta + 1)$-vertex-coloring or $(2\Delta-1)$-edge-coloring algorithm that runs in $\mathrm{poly}(\log n)$ rounds regardless of the value of $\Delta$ has been a major challenge \cite[Problems 11.2 and 11.4]{BE}. After a long line of contributions by numerous researchers, such algorithms have been discovered by Fischer, Ghaffari, and Kuhn \cite{FGK} for edge-coloring and by Rozho\v{n} and Ghaffari \cite{RG} for vertex-coloring.
	
	There have also been many results concerning edge-coloring with fewer than $2\Delta - 1$ colors. For thorough surveys see, e.g., \cite{CHLPU,GKMU}. Here we only mention a few highlights. Chang, He, Li, Pettie, and Uitto \cite{CHLPU} developed a randomized algorithm for $(\Delta + \tilde{O}(\sqrt{\Delta}))$-edge-coloring that runs in $\mathrm{poly}(\Delta, \log \log n)$ rounds. \ep{Their algorithm relies on efficiently solving instances of the Lov\'asz Local Lemma, which can be done in $\mathrm{poly}(\Delta, \log \log n)$ due to a recent result of Rozho\v{n} and Ghaffari \cite{RG}.} In the deterministic setting, Ghaffari, Kuhn, Maus, and Uitto \cite{GKMU} gave a $\mathrm{poly}(\Delta, \log n)$-round algorithm for $\lfloor 3\Delta/2\rfloor$-edge-coloring. The smallest number of colors for which an efficient distributed algorithm was known prior to our work is $\Delta + 2$ \ep{i.e., just one more than the bound in Vizing's theorem}. Namely, a randomized $\mathrm{poly}(\Delta, \log n)$-round algorithm for $(\Delta+2)$\-/edge\-/coloring was designed by Su and Vu \cite{SV}. In a recent breakthrough, Rozho\v{n} and Ghaffari \cite{RG} developed a general derandomization technique that, in particular, allows one to modify the Su--Vu algorithm to make it deterministic.
	
	It is natural to wonder how close to optimal the number of rounds required by our algorithm is. Chang, He, Li, Pettie, and Uitto \cite{CHLPU} showed that any $(\Delta+1)$-edge-coloring algorithm based on ``extending partial colorings by recoloring subgraphs'' \ep{a class to which our algorithm belongs} must take at least $\Omega(\Delta \log (n/\Delta))$ rounds.	Furthermore, they showed that every deterministic \LOCAL algorithm for $(2\Delta-2)$-edge-coloring requires $\Omega(\log n/\log \Delta)$ rounds, even in the case when the underlying graph $G$ is a tree. The number of rounds in our algorithm is certainly not linear in $\Delta \log n$, although we did not make an attempt to optimize it precisely. The following remains an interesting open problem:
	
	\begin{ques}
		In the regime when $\Delta$ is constant, does there exist a randomized $(\Delta+1)$-edge-coloring algorithm that takes only $o(\log n)$ rounds?
	\end{ques}
	
	\subsection{Locality of $(\Delta+1)$-edge-coloring}\label{subsec:intro_comb}
	
	The main ingredient of our algorithm is a certain purely combinatorial result concerning extensions of partial $(\Delta+1)$-edge-colorings, namely Theorem~\ref{theo:ext} below. Before stating it, we need to introduce some terminology and explain our motivation.
	Let $G$ be an $n$-vertex graph of maximum degree $\Delta$. Fix a partition $E(G) = X \sqcup U$ of the edge set of $G$ and let $\phi \colon X \to [\Delta + 1]$ be a proper partial edge-coloring with domain $X$. We call the edges in $X$ \ep{resp.~$U$} \emphd{colored} \ep{resp.~\emphd{uncolored}} by $\phi$. A subgraph $H \subseteq G$ is \emphd{augmenting} \ep{for $\phi$} if at least one edge of $H$ is uncolored and there is a proper coloring $\psi \colon X \cup E(H) \to [\Delta + 1]$ that agrees with $\phi$ on $X \setminus E(H)$; that is, by only modifying $\phi$ on the edges of $H$, it is possible to extend it to a proper partial coloring with domain $X \cup E(H)$. Notice that if $H_1$, \ldots, $H_k$ are vertex-disjoint augmenting subgraphs, then their \ep{disjoint} union $H_1 \sqcup \ldots \sqcup H_k$ is also augmenting.
	
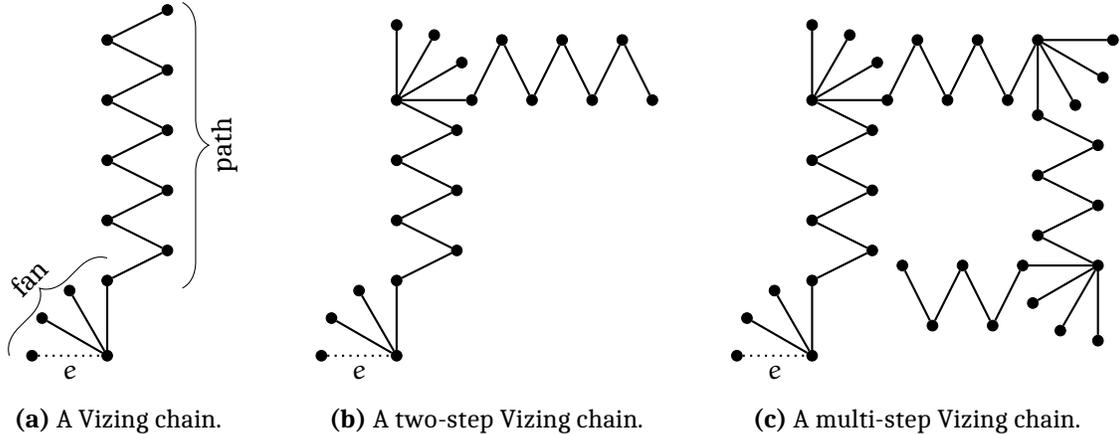
\begin{figure}[t]
		\centering
		\begin{subfigure}[t]{.2\textwidth}
			\centering
			\begin{tikzpicture}
			\node[circle,fill=black,draw,inner sep=0pt,minimum size=4pt] (a) at (0,0) {};
			\node[circle,fill=black,draw,inner sep=0pt,minimum size=4pt] (b) at (1,0) {};
			\path (b) ++(150:1) node[circle,fill=black,draw,inner sep=0pt,minimum size=4pt] (c) {};
			\path (b) ++(120:1) node[circle,fill=black,draw,inner sep=0pt,minimum size=4pt] (d) {};
			\node[circle,fill=black,draw,inner sep=0pt,minimum size=4pt] (e) at (1,1) {};
			\node[circle,fill=black,draw,inner sep=0pt,minimum size=4pt] (f) at (1.8,1.4) {};
			\node[circle,fill=black,draw,inner sep=0pt,minimum size=4pt] (g) at (1,1.8) {};
			\node[circle,fill=black,draw,inner sep=0pt,minimum size=4pt] (h) at (1.8,2.2) {};
			\node[circle,fill=black,draw,inner sep=0pt,minimum size=4pt] (i) at (1,2.6) {};
			\node[circle,fill=black,draw,inner sep=0pt,minimum size=4pt] (j) at (1.8,3) {};
			\node[circle,fill=black,draw,inner sep=0pt,minimum size=4pt] (k) at (1,3.4) {};
			\node[circle,fill=black,draw,inner sep=0pt,minimum size=4pt] (l) at (1.8,3.8) {};
			\node[circle,fill=black,draw,inner sep=0pt,minimum size=4pt] (m) at (1,4.2) {};
			\node[circle,fill=black,draw,inner sep=0pt,minimum size=4pt] (n) at (1.8,4.6) {};
			
			\draw[ thick,dotted] (a) to node[midway,anchor=north] {$e$} (b);
			\draw[ thick] (b) -- (c) (b) -- (d) (b) -- (e) -- (f) -- (g) -- (h) -- (i) -- (j) -- (k) -- (l) -- (m) -- (n);
			
			\draw[decoration={brace,amplitude=10pt,mirror},decorate] (2,0.9) -- node [midway,below,sloped,yshift=-7pt] {path} (2,4.7);
			
			\draw[decoration={brace,amplitude=10pt},decorate] (-0.3,0) -- node [midway,above,sloped,yshift=8pt] {fan} (1,1.3);
			\end{tikzpicture}
			\caption{A Vizing chain.}\label{fig:chains:Vizing}
		\end{subfigure}%
		\qquad%
		\begin{subfigure}[t]{.3\textwidth}
			\centering
			\begin{tikzpicture}
			\node[circle,fill=black,draw,inner sep=0pt,minimum size=4pt] (a) at (0,0) {};
			\node[circle,fill=black,draw,inner sep=0pt,minimum size=4pt] (b) at (1,0) {};
			\path (b) ++(150:1) node[circle,fill=black,draw,inner sep=0pt,minimum size=4pt] (c) {};
			\path (b) ++(120:1) node[circle,fill=black,draw,inner sep=0pt,minimum size=4pt] (d) {};
			\node[circle,fill=black,draw,inner sep=0pt,minimum size=4pt] (e) at (1,1) {};
			\node[circle,fill=black,draw,inner sep=0pt,minimum size=4pt] (f) at (1.8,1.4) {};
			\node[circle,fill=black,draw,inner sep=0pt,minimum size=4pt] (g) at (1,1.8) {};
			\node[circle,fill=black,draw,inner sep=0pt,minimum size=4pt] (h) at (1.8,2.2) {};
			\node[circle,fill=black,draw,inner sep=0pt,minimum size=4pt] (i) at (1,2.6) {};
			\node[circle,fill=black,draw,inner sep=0pt,minimum size=4pt] (j) at (1.8,3) {};
			\node[circle,fill=black,draw,inner sep=0pt,minimum size=4pt] (k) at (1,3.4) {};
			\node[circle,fill=black,draw,inner sep=0pt,minimum size=4pt] (l) at (1,4.4) {};
			
			
			\path (k) ++(60:1) node[circle,fill=black,draw,inner sep=0pt,minimum size=4pt] (m) {};
			\path (k) ++(30:1) node[circle,fill=black,draw,inner sep=0pt,minimum size=4pt] (n) {};
			\path (k) ++(0:1) node[circle,fill=black,draw,inner sep=0pt,minimum size=4pt] (o) {};
			
			\node[circle,fill=black,draw,inner sep=0pt,minimum size=4pt] (p) at (2.4,4.2) {};
			\node[circle,fill=black,draw,inner sep=0pt,minimum size=4pt] (q) at (2.8,3.4) {};
			\node[circle,fill=black,draw,inner sep=0pt,minimum size=4pt] (r) at (3.2,4.2) {};
			\node[circle,fill=black,draw,inner sep=0pt,minimum size=4pt] (s) at (3.6,3.4) {};
			\node[circle,fill=black,draw,inner sep=0pt,minimum size=4pt] (t) at (4,4.2) {};
			\node[circle,fill=black,draw,inner sep=0pt,minimum size=4pt] (u) at (4.4,3.4) {};
			
			\draw[ thick,dotted] (a) to node[midway,anchor=north] {$e$} (b);
			\draw[ thick] (b) -- (c) (b) -- (d) (b) -- (e) -- (f) -- (g) -- (h) -- (i) -- (j) -- (k) -- (l);
			
			\draw[ thick] (k) -- (m) (k) -- (n) (k) -- (o) -- (p) -- (q) -- (r) -- (s) -- (t) -- (u);
			\end{tikzpicture}
			\caption{A two-step Vizing chain.}\label{fig:chains:two}
		\end{subfigure}%
		\qquad%
		\begin{subfigure}[t]{.3\textwidth}
			\centering
			\begin{tikzpicture}
			\node[circle,fill=black,draw,inner sep=0pt,minimum size=4pt] (a) at (0,0) {};
			\node[circle,fill=black,draw,inner sep=0pt,minimum size=4pt] (b) at (1,0) {};
			\path (b) ++(150:1) node[circle,fill=black,draw,inner sep=0pt,minimum size=4pt] (c) {};
			\path (b) ++(120:1) node[circle,fill=black,draw,inner sep=0pt,minimum size=4pt] (d) {};
			\node[circle,fill=black,draw,inner sep=0pt,minimum size=4pt] (e) at (1,1) {};
			\node[circle,fill=black,draw,inner sep=0pt,minimum size=4pt] (f) at (1.8,1.4) {};
			\node[circle,fill=black,draw,inner sep=0pt,minimum size=4pt] (g) at (1,1.8) {};
			\node[circle,fill=black,draw,inner sep=0pt,minimum size=4pt] (h) at (1.8,2.2) {};
			\node[circle,fill=black,draw,inner sep=0pt,minimum size=4pt] (i) at (1,2.6) {};
			\node[circle,fill=black,draw,inner sep=0pt,minimum size=4pt] (j) at (1.8,3) {};
			\node[circle,fill=black,draw,inner sep=0pt,minimum size=4pt] (k) at (1,3.4) {};
			\node[circle,fill=black,draw,inner sep=0pt,minimum size=4pt] (l) at (1,4.4) {};
			
			\path (k) ++(60:1) node[circle,fill=black,draw,inner sep=0pt,minimum size=4pt] (m) {};
			\path (k) ++(30:1) node[circle,fill=black,draw,inner sep=0pt,minimum size=4pt] (n) {};
			\path (k) ++(0:1) node[circle,fill=black,draw,inner sep=0pt,minimum size=4pt] (o) {};
			
			\node[circle,fill=black,draw,inner sep=0pt,minimum size=4pt] (p) at (2.4,4.2) {};
			\node[circle,fill=black,draw,inner sep=0pt,minimum size=4pt] (q) at (2.8,3.4) {};
			\node[circle,fill=black,draw,inner sep=0pt,minimum size=4pt] (r) at (3.2,4.2) {};
			\node[circle,fill=black,draw,inner sep=0pt,minimum size=4pt] (s) at (3.6,3.4) {};
			\node[circle,fill=black,draw,inner sep=0pt,minimum size=4pt] (t) at (4,4.2) {};

			\node[circle,fill=black,draw,inner sep=0pt,minimum size=4pt] (v) at (5,4.2) {};
			\path (t) ++(-30:1) node[circle,fill=black,draw,inner sep=0pt,minimum size=4pt] (w) {};
			\path (t) ++(-60:1) node[circle,fill=black,draw,inner sep=0pt,minimum size=4pt] (x) {};
			\path (t) ++(-90:1) node[circle,fill=black,draw,inner sep=0pt,minimum size=4pt] (y) {};

			\node[circle,fill=black,draw,inner sep=0pt,minimum size=4pt] (z) at (4.8,2.8) {};
			\node[circle,fill=black,draw,inner sep=0pt,minimum size=4pt] (aa) at (4,2.4) {};
			\node[circle,fill=black,draw,inner sep=0pt,minimum size=4pt] (ab) at (4.8,2) {};
			\node[circle,fill=black,draw,inner sep=0pt,minimum size=4pt] (ac) at (4,1.6) {};
			\node[circle,fill=black,draw,inner sep=0pt,minimum size=4pt] (ad) at (4.8,1.2) {};
			\node[circle,fill=black,draw,inner sep=0pt,minimum size=4pt] (ae) at (4.8,0.2) {};
			
			\path (ad) ++(-120:1) node[circle,fill=black,draw,inner sep=0pt,minimum size=4pt] (ag) {};
			\path (ad) ++(-150:1) node[circle,fill=black,draw,inner sep=0pt,minimum size=4pt] (ah) {};
			\path (ad) ++(-180:1) node[circle,fill=black,draw,inner sep=0pt,minimum size=4pt] (ai) {};
			
			\node[circle,fill=black,draw,inner sep=0pt,minimum size=4pt] (aj) at (3.4,0.4) {};
			\node[circle,fill=black,draw,inner sep=0pt,minimum size=4pt] (ak) at (3,1.2) {};
			\node[circle,fill=black,draw,inner sep=0pt,minimum size=4pt] (al) at (2.6,0.4) {};
			\node[circle,fill=black,draw,inner sep=0pt,minimum size=4pt] (am) at (2.2,1.2) {};

			\draw[ thick,dotted] (a) to node[midway,anchor=north] {$e$} (b);
			\draw[ thick] (b) -- (c) (b) -- (d) (b) -- (e) -- (f) -- (g) -- (h) -- (i) -- (j) -- (k) -- (l);
			
			\draw[ thick] (k) -- (m) (k) -- (n) (k) -- (o) -- (p) -- (q) -- (r) -- (s) -- (t) -- (v);
			
			\draw[ thick] (t) -- (w) (t) -- (x) (t) -- (y) -- (z) -- (aa) -- (ab) -- (ac) -- (ad) -- (ae);
			
			\draw[ thick] (ad) -- (ag) (ad) -- (ah) (ad) -- (ai) -- (aj) -- (ak) -- (al) -- (am);
			\end{tikzpicture}
			\caption{A multi-step Vizing chain.}\label{fig:chains:multi}
		\end{subfigure}%
		\caption{Types of augmenting subgraphs.}\label{fig:chains}
	\end{figure}

	A standard way to prove Vizing's theorem is to construct, given an uncolored edge $e \in U$, an augmenting graph $H$ with $e \in E(H)$ of a certain special form, called a \emph{Vizing chain} \ep{see Fig.~\ref{fig:chains:Vizing}}. A Vizing chain $H$ consists of a \emph{fan}---i.e., a set of edges that share the same common vertex---and a \emph{path} that is \emph{alternating} in the sense that the sequence of colors assigned by $\phi$ to its edges has the form $\alpha$, $\beta$, $\alpha$, $\beta$, \ldots{} for a pair of colors $\alpha$, $\beta$. In the distributed setting, this approach raises two difficulties. First, uncolored edges cannot be treated one at a time; it is necessary to be able to extend the coloring to include a non-negligible fraction of the uncolored edges simultaneously. The second difficulty is even more fundamental: the Vizing chain $H$ corresponding to an uncolored edge $e$ may have large diameter \ep{in principle, linear in $n$}, so it cannot be discovered in a small number of rounds in the \LOCAL model. We overcome both these difficulties by showing that it is possible to find a ``large'' collection of vertex-disjoint ``small'' augmenting subgraphs \ep{although these subgraphs need not be Vizing chains anymore}:
	
	\begin{theo}\label{theo:ext}
		There is $n_0 \in \N$ such that for all $n > n_0$, the following holds. Let $G$ be an $n$-vertex graph of maximum degree $\Delta$. Fix a partition $E(G) = X \sqcup U$ and let $\phi \colon X \to [\Delta + 1]$ be a proper partial edge-coloring. Then there exists a subset $W \subseteq U$ of size $|W| \geq |U|/((\Delta+1)^{10}(\log n)^2)$ such that it is possible to assign to each edge $e \in W$ a connected augmenting subgraph $H^e \subseteq G$ with the following properties:
		\begin{itemize}
			\item for each $e \in W$, 
			$E(H^e) \cap U = \set{e}$, that is, $e$ is the unique uncolored edge in $H^e$;
			\item the vertex sets of the graphs $H^e$, $e \in W$, are pairwise disjoint;
			\item for each $e \in W$, $|E(H^e)| \leq (\Delta+1)^6(\log n)^2$.
		\end{itemize}
	\end{theo}
	
	Thanks to the recent work of Ghaffari, Harris, and Kuhn \cite{derand} and Harris \cite{Harris} on distributed approximation algorithms for hypergraph maximum matching, it is straightforward to derive Theorem~\ref{theo:main} from Theorem~\ref{theo:ext}. We present the details of this derivation in \S\ref{sec:Harris}. Note that if $U \neq \0$, then the lower bound on the size of $W$ in Theorem~\ref{theo:ext} is positive. Thus, Theorem~\ref{theo:ext} in particular implies that if $G$ contains only a single uncolored edge $e$, then there is a connected augmenting subgraph $H \subseteq G$ such that $e \in E(H)$ and $|E(H)| \leq \mathrm{poly}(\Delta, \log n)$.
	
	Let us now say a few words about the proof of Theorem~\ref{theo:ext}. Our inspiration comes from recent developments in the area of \emph{descriptive combinatorics}, i.e., the study of combinatorial problems on infinite graphs in the presence of additional topological or measure-theoretic constraints. For a state-of-the-art introduction to this subject, see the surveys \cite{KechrisMarks} by Kechris and Marks and \cite{Pikh_survey} by Pikhurko. It turns out that many of the challenges one encounters in descriptive combinatorics mirror those in distributed computing; furthermore, there are some explicit implications between the two areas \cite{BerDist, BerCont}. Until very recently, one of the most celebrated open problems in descriptive combinatorics had been to obtain a ``measurable'' version of Vizing's theorem for graphs $G$ on standard probability spaces \ep{for instance, with $V(G) = [0,1]$, the unit interval}. This was recently accomplished \ep{under some minor additional assumptions} by Greb\'\i{}k and Pikhurko \cite{GP}. Their key idea was to consider a more general type of augmenting graphs, namely \emph{two-step Vizing chains} \ep{called ``iterated Vizing chains'' in \cite{GP}; see Fig.~\ref{fig:chains:two}}. Roughly speaking, to construct a two-step Vizing chain, one first starts growing a Vizing chain as usual, but then at some point, instead of completing the alternating path \ep{which may be too long}, one ``changes the direction'' and starts growing a second Vizing chain. The flexibility in the choice of the point where the second Vizing chain starts allowed Greb\'\i{}k and Pikhurko to control the total number of edges in the resulting structure. In particular, returning to finite graphs, a relatively straightforward adaptation of the Greb\'\i{}k--Pikhurko argument shows that for each uncolored edge $e \in U$, one can find an augmenting two-step Vizing chain $H$ with $e \in E(H)$ such that $|E(H)| = O(\mathrm{poly}(\Delta)\sqrt{n})$.
	
	Unfortunately, the bound $O(\mathrm{poly}(\Delta)\sqrt{n})$ is not sufficient for our purposes, as we can only work with augmenting graphs of diameter at most $\mathrm{poly}(\Delta, \log n)$. We solve this issue by extending the ideas of Greb\'\i{}k and Pikhurko and considering \emph{multi-step Vizing chains} \ep{see Fig.~\ref{fig:chains:multi}}. Roughly speaking, we prove Theorem~\ref{theo:ext} by showing that it is possible to construct augmenting subgraphs of this form by joining together $O(\log n)$ partial Vizing chains, each consisting of $O(\mathrm{poly}(\Delta)\log n)$ edges. For the reader's convenience, we give an informal overview of our argument in \S\ref{sec:overview}.
	
	\subsubsection*{{Acknowledgments}}
	
	I am very grateful to Clinton Conley and Oleg Pikhurko for insightful conversations and to the anonymous referees for carefully reading the manuscript and providing helpful suggestions.
	
	\section{From Theorem~\ref{theo:ext} to Theorem~\ref{theo:main}}\label{sec:Harris}
	
	In this section we deduce Theorem~\ref{theo:main} from Theorem~\ref{theo:ext}. To this end, we shall use the distributed approximation algorithm for hypergraph maximum matching due to Harris \cite{Harris}. \ep{Alternatively, we could have used the earlier and somewhat less efficient algorithm due to Ghaffari, Harris, and Kuhn \cite[Theorem 1.2]{derand}.} To begin with, we need a few definitions and some notation. Recall that $\log^\ast n$ denotes the iterated logarithm of $n$. 
	The asymptotic notation $\tilde{O}(x)$ hides polylogarithmic factors, i.e., it stands for $O( \mathrm{poly}(\log x) \cdot x)$. Let $\mathscr{H} = (V(\mathscr{H}), E(\mathscr{H}))$ be a hypergraph. To avoid potential confusion with the graph case, we will call the elements of $E(\mathscr{H})$ \emphd{hyperedges}. The \emphd{rank} of $\mathscr{H}$, in symbols $r(\mathscr{H})$, is the largest size of a hyperedge of $\mathscr{H}$ \ep{if $E(\mathscr{H}) = \0$, then we set $r(\mathscr{H}) \defeq 0$}. The \emphd{maximum degree} $d(\mathscr{H})$ of $\mathscr{H}$ is the maximum, over all the vertices $x \in V(\mathscr{H})$, of the number of hyperedges $S \in E(\mathscr{H})$ with $x \in S$. A \emphd{matching} in $\mathscr{H}$ is a set $M \subseteq E(\mathscr{H})$ of pairwise disjoint hyperedges. We use $\mu(\mathscr{H})$ to denote the maximum number of edges in a matching in $\mathscr{H}$. While the usual \LOCAL model is defined for graphs, there is an analogous model operating on a hypergraph $\mathscr{H}$. Namely, in a single communication round of the \LOCAL model on $\mathscr{H}$, each vertex $x \in V(\mathscr{H})$ is allowed to send messages to every vertex $y \in V(\mathscr{H})$ such that $x$ and $y$ are contained in a common hyperedge. 
	
	\begin{theo}[{Harris \cite[Theorem 1.1]{Harris}}]\label{theo:Harris}
		There exists a deterministic distributed algorithm in the \LOCAL model on an $n$-vertex hypergraph $\mathscr{H}$ that outputs a matching $M \subseteq E(\mathscr{H})$ with $|M| = \Omega(\mu/r)$ in \[\tilde{O}(r \log d + \log^2 d + \log^\ast n)\] rounds, where $r \defeq r(\mathscr{H})$, $d \defeq d(\mathscr{H})$, and $\mu \defeq \mu(\mathscr{H})$.
	\end{theo}
	
	With Theorem~\ref{theo:Harris} in hand, we are ready to derive Theorem~\ref{theo:main} from Theorem~\ref{theo:ext}. Let $G$ be an $n$\=/vertex graph of maximum degree $\Delta$. Since we only care about asymptotic bounds, we may assume that $n$ is sufficiently large for Theorem~\ref{theo:ext} to apply. For brevity, set $V \defeq V(G)$ and $E \defeq E(G)$. We shall devise a $\mathrm{poly}(\Delta, \log n)$-round \LOCAL algorithm that, given a proper partial edge-coloring $\phi \colon X \to [\Delta + 1]$, $X \subseteq E$, computes a new proper partial edge-coloring $\psi \colon Y \to [\Delta + 1]$, $Y \subseteq E$, such that $Y \supseteq X$ and
	\[
	|E \setminus Y| \,\leq\, (1 - 1/\mathrm{poly}(\Delta, \log n)) |E \setminus X|.
	\]
	Iterating this algorithm $\mathrm{poly}(\Delta, \log n)$ times starting with the empty coloring will then yield a proper $(\Delta + 1)$-edge-coloring of the entire graph $G$, as desired.
	
	So, let us fix a proper partial edge-coloring $\phi \colon X \to [\Delta + 1]$. Define an auxiliary hypergraph $\mathscr{H}$ on the vertex set $V$ by making a subset $S \subseteq V$ a hyperedge of $\mathscr{H}$ if and only if $S = V(H)$ for some connected augmenting subgraph $H \subseteq G$ with $|E(H)| \leq (\Delta+1)^6 (\log n)^2$. Theorem~\ref{theo:ext} then implies that
	\[
	\mu(\mathscr{H}) \,\geq\, |E \setminus X|/((\Delta+1)^{10} (\log n)^2) \,=\, |E \setminus X|/\mathrm{poly}(\Delta, \log n).
	\]
	By definition, each $S \in E(\mathscr{H})$ satisfies $|S| \leq (\Delta+1)^6 (\log n)^2 + 1 = \mathrm{poly}(\Delta, \log n)$, so
	\[
	r(\mathscr{H}) \,\leq\, \mathrm{poly}(\Delta, \log n).
	\]
	To bound the maximum degree of $\mathscr{H}$, consider any vertex $x \in V$. If $x \in S \in E(\mathscr{H})$, then the vertices of $S$ can be ordered as $x = x_0$, $x_1$, \ldots, $x_{|S|-1}$ so that each $x_i$, $i \geq 1$, is adjacent to at least one of $x_0$, \ldots, $x_{i-1}$. This means that once $x_0$, \ldots, $x_{i-1}$ are fixed, there are at most $i \Delta \leq r(\mathscr{H})\Delta$ choices for $x_i$, and thus
	\[
	d(\mathscr{H}) \,\leq\, (r(\mathscr{H})\Delta)^{r(\mathscr{H})} \,\leq\, \exp(\mathrm{poly}(\Delta, \log n)).
	\]
	Therefore, Theorem~\ref{theo:Harris} provides a \LOCAL algorithm on $\mathscr{H}$ that outputs a matching $M \subseteq E(\mathscr{H})$ with
	\[
	|M| \,\geq\, |E \setminus X|/ \mathrm{poly}(\Delta, \log n)
	\]
	in $\mathrm{poly}(\Delta, \log n)$ rounds. It is clear that a single round of communication in the \LOCAL model on $\mathscr{H}$ can be simulated in $\mathrm{poly}(\Delta, \log n)$ rounds of the \LOCAL model on $G$, so $M$ can also be computed in the \LOCAL model on $G$ in $\mathrm{poly}(\Delta, \log n)$ rounds. Once such a matching $M$ is fixed, in $\mathrm{poly}(\Delta, \log n)$ rounds it is possible to choose a ``leader'' in each hyperedge $S \in M$ \ep{by picking, say, the vertex with the lexicographically least identifier}. In another $\mathrm{poly}(\Delta, \log n)$ rounds, each ``leader'' $x$ surveys its corresponding hyperedge $S$, picks an arbitrary augmenting subgraph $H$ with $V(H) = S$, and modifies the coloring so that every edge of $H$ becomes colored. Since $M$ is a matching in $\mathscr{H}$, all these modifications can happen in parallel without conflicting with each other. In the resulting coloring, every edge in $X$ is still colored, plus there is at least one new colored edge per each $S \in M$, and hence we are done.

	\section{Informal overview of the proof of Theorem~\ref{theo:ext}}\label{sec:overview}
	
	In this section we sketch the main ideas behind our proof of Theorem \ref{theo:ext}. Due to the informal character of this section, we will avoid any technicalities and precise definitions, all of which are postponed to the later sections. It should be understood that the terminology we use here might have a slightly different meaning when we formally define it, but only in minor and technical ways.
	
	Let $G$ be an $n$-vertex graph of maximum degree $\Delta$ and let $\phi$ be a proper partial $(\Delta+1)$-edge-coloring of $G$. For simplicity, in this section we shall treat $\Delta$ as a constant, so that implicit constants in asymptotic notation may depend on $\Delta$. Fix an uncolored edge $e = xy$. We will describe our strategy for proving the following consequence of Theorem~\ref{theo:ext}:
	\[
	\parbox{0.95\textwidth}{\emph{There is a connected augmenting subgraph $H \subseteq G$ such that $e \in E(H)$ and $|E(H)| \leq O((\log n)^2)$.}}
	\]
	
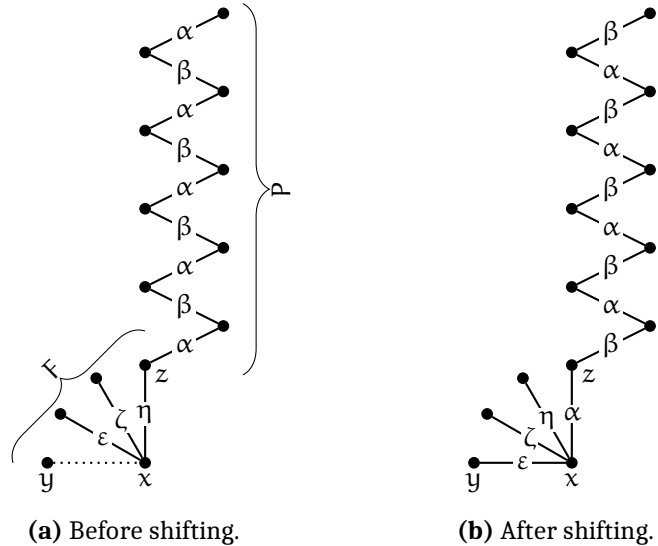
\begin{figure}[h]
		\centering
		\begin{subfigure}[t]{.2\textwidth}
			\centering
			\begin{tikzpicture}[scale=1.3]
			\node[circle,fill=black,draw,inner sep=0pt,minimum size=4pt] (a) at (0,0) {};
			\node[circle,fill=black,draw,inner sep=0pt,minimum size=4pt] (b) at (1,0) {};
			\path (b) ++(150:1) node[circle,fill=black,draw,inner sep=0pt,minimum size=4pt] (c) {}; 
			\path (b) ++(120:1) node[circle,fill=black,draw,inner sep=0pt,minimum size=4pt] (d) {}; 
			\node[circle,fill=black,draw,inner sep=0pt,minimum size=4pt] (e) at (1,1) {};
			\node[circle,fill=black,draw,inner sep=0pt,minimum size=4pt] (f) at (1.8,1.4) {};
			\node[circle,fill=black,draw,inner sep=0pt,minimum size=4pt] (g) at (1,1.8) {};
			\node[circle,fill=black,draw,inner sep=0pt,minimum size=4pt] (h) at (1.8,2.2) {};
			\node[circle,fill=black,draw,inner sep=0pt,minimum size=4pt] (i) at (1,2.6) {};
			\node[circle,fill=black,draw,inner sep=0pt,minimum size=4pt] (j) at (1.8,3) {};
			\node[circle,fill=black,draw,inner sep=0pt,minimum size=4pt] (k) at (1,3.4) {};
			\node[circle,fill=black,draw,inner sep=0pt,minimum size=4pt] (l) at (1.8,3.8) {};
			\node[circle,fill=black,draw,inner sep=0pt,minimum size=4pt] (m) at (1,4.2) {};
			\node[circle,fill=black,draw,inner sep=0pt,minimum size=4pt] (n) at (1.8,4.6) {};
			
			\draw[ thick,dotted] (a) -- (b);
			\draw[ thick] (b) to node[midway,inner sep=1pt,outer sep=1pt,minimum size=4pt,fill=white] {$\epsilon$} (c) (b) to node[midway,inner sep=0.5pt,outer sep=0.5pt,minimum size=4pt,fill=white] {$\zeta$} (d) (b) to node[midway,inner sep=1pt,outer sep=1pt,minimum size=4pt,fill=white] {$\eta$} (e) to node[midway,inner sep=1pt,outer sep=1pt,minimum size=4pt,fill=white] {$\alpha$} (f) to node[midway,inner sep=1pt,outer sep=1pt,minimum size=4pt,fill=white] {$\beta$} (g) to node[midway,inner sep=1pt,outer sep=1pt,minimum size=4pt,fill=white] {$\alpha$} (h) to node[midway,inner sep=1pt,outer sep=1pt,minimum size=4pt,fill=white] {$\beta$} (i) to node[midway,inner sep=1pt,outer sep=1pt,minimum size=4pt,fill=white] {$\alpha$} (j) to node[midway,inner sep=1pt,outer sep=1pt,minimum size=4pt,fill=white] {$\beta$} (k) to node[midway,inner sep=1pt,outer sep=1pt,minimum size=4pt,fill=white] {$\alpha$} (l) to node[midway,inner sep=1pt,outer sep=1pt,minimum size=4pt,fill=white] {$\beta$} (m) to node[midway,inner sep=1pt,outer sep=1pt,minimum size=4pt,fill=white] {$\alpha$} (n);
			
			\draw[decoration={brace,amplitude=10pt,mirror},decorate] (2,0.9) -- node [midway,below,sloped,yshift=-7pt] {$P$} (2,4.7);
			
			\draw[decoration={brace,amplitude=10pt},decorate] (-0.35,0) -- node [midway,above,sloped,yshift=7pt] {$F$} (1,1.35);
			
			\node[anchor=148.28] at (e) {$z$};
			\node[anchor=north] at (a) {$y$};
			\node[anchor=north] at (b) {$x$};
			\end{tikzpicture}
			\caption{Before shifting.}\label{fig:Viz_col:initial}
		\end{subfigure}%
		\qquad\qquad\qquad%
		\begin{subfigure}[t]{.2\textwidth}
			\centering
			\begin{tikzpicture}[scale=1.3]
			\node[circle,fill=black,draw,inner sep=0pt,minimum size=4pt] (a) at (0,0) {};
			\node[circle,fill=black,draw,inner sep=0pt,minimum size=4pt] (b) at (1,0) {};
			\path (b) ++(150:1) node[circle,fill=black,draw,inner sep=0pt,minimum size=4pt] (c) {}; 
			\path (b) ++(120:1) node[circle,fill=black,draw,inner sep=0pt,minimum size=4pt] (d) {}; 
			\node[circle,fill=black,draw,inner sep=0pt,minimum size=4pt] (e) at (1,1) {};
			\node[circle,fill=black,draw,inner sep=0pt,minimum size=4pt] (f) at (1.8,1.4) {};
			\node[circle,fill=black,draw,inner sep=0pt,minimum size=4pt] (g) at (1,1.8) {};
			\node[circle,fill=black,draw,inner sep=0pt,minimum size=4pt] (h) at (1.8,2.2) {};
			\node[circle,fill=black,draw,inner sep=0pt,minimum size=4pt] (i) at (1,2.6) {};
			\node[circle,fill=black,draw,inner sep=0pt,minimum size=4pt] (j) at (1.8,3) {};
			\node[circle,fill=black,draw,inner sep=0pt,minimum size=4pt] (k) at (1,3.4) {};
			\node[circle,fill=black,draw,inner sep=0pt,minimum size=4pt] (l) at (1.8,3.8) {};
			\node[circle,fill=black,draw,inner sep=0pt,minimum size=4pt] (m) at (1,4.2) {};
			\node[circle,fill=black,draw,inner sep=0pt,minimum size=4pt] (n) at (1.8,4.6) {};
			
			\draw[ thick] (a) to node[midway,inner sep=1pt,outer sep=1pt,minimum size=4pt,fill=white] {$\epsilon$} (b);
			\draw[ thick] (b) to node[midway,inner sep=0.5pt,outer sep=0.5pt,minimum size=4pt,fill=white] {$\zeta$} (c) (b) to node[midway,inner sep=1pt,outer sep=1pt,minimum size=4pt,fill=white] {$\eta$} (d) (b) to node[midway,inner sep=1pt,outer sep=1pt,minimum size=4pt,fill=white] {$\alpha$} (e) to node[midway,inner sep=1pt,outer sep=1pt,minimum size=4pt,fill=white] {$\beta$} (f) to node[midway,inner sep=1pt,outer sep=1pt,minimum size=4pt,fill=white] {$\alpha$} (g) to node[midway,inner sep=1pt,outer sep=1pt,minimum size=4pt,fill=white] {$\beta$} (h) to node[midway,inner sep=1pt,outer sep=1pt,minimum size=4pt,fill=white] {$\alpha$} (i) to node[midway,inner sep=1pt,outer sep=1pt,minimum size=4pt,fill=white] {$\beta$} (j) to node[midway,inner sep=1pt,outer sep=1pt,minimum size=4pt,fill=white] {$\alpha$} (k) to node[midway,inner sep=1pt,outer sep=1pt,minimum size=4pt,fill=white] {$\beta$} (l) to node[midway,inner sep=1pt,outer sep=1pt,minimum size=4pt,fill=white] {$\alpha$} (m) to node[midway,inner sep=1pt,outer sep=1pt,minimum size=4pt,fill=white] {$\beta$} (n);
			
			
			
			\node[anchor=148.28] at (e) {$z$};
			\node[anchor=north] at (a) {$y$};
			\node[anchor=north] at (b) {$x$};
			\end{tikzpicture}
			\caption{After shifting.}\label{fig:Viz_col:shifted}
		\end{subfigure}%
		\caption{A Vizing chain $H = F + P$ before and after shifting. Greek letters represent colors.}\label{fig:Viz_col}
	\end{figure}

	As mentioned in \S\ref{subsec:intro_comb}, standard proofs of Vizing's theorem \ep{for instance, the one given in \cite[\S5.3]{Die}} proceed by building an augmenting subgraph $H$ consisting of two parts:
	\begin{itemize}
		\item a \emph{fan} $F$ with pivot $x$, i.e., a sequence of distinct edges incident to $x$, starting with $e$;
		
		\item an \emph{$\alpha\beta$-alternating path} $P$, i.e., a path whose edges are colored $\alpha$, $\beta$, $\alpha$, $\beta$, etc.
	\end{itemize} 
	There is some variation between different presentations of the proof of Vizing's theorem in the literature as to the precise way in which the fan $F$ and the path $P$ fit together. For our purposes, it will be most convenient to assume that $P$ starts at the vertex $z$, where $xz$ is the last edge of $F$. This situation is illustrated in Fig.~\ref{fig:Viz_col:initial}. We write $H = F + P$ to indicate that the graph $H$ is built from $F$ and $P$ in this manner and refer to $H$ as a \emph{Vizing chain}. Given a Vizing chain $H = F + P$, we can \emph{shift} the colors in $H$, as shown in Fig.~\ref{fig:Viz_col:shifted}. One can show \ep{see, e.g., \cite[\S5.3]{Die} or Corollary~\ref{corl:vizing} below} that there is such a Vizing chain $H$ that is \emph{shiftable}, meaning that the coloring resulting after the shifting is proper, and thus $H$ is an augmenting subgraph.
	
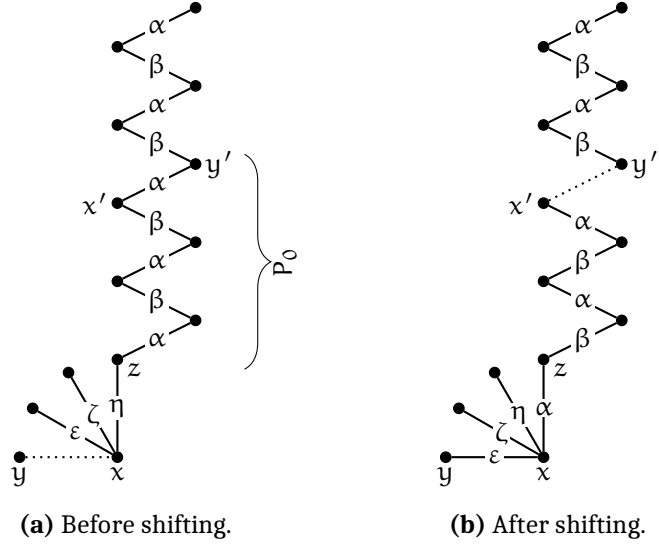
\begin{figure}[h]
		\centering
		\begin{subfigure}[t]{.2\textwidth}
			\centering
			\begin{tikzpicture}[scale=1.3]
			\node[circle,fill=black,draw,inner sep=0pt,minimum size=4pt] (a) at (0,0) {};
			\node[circle,fill=black,draw,inner sep=0pt,minimum size=4pt] (b) at (1,0) {};
			\path (b) ++(150:1) node[circle,fill=black,draw,inner sep=0pt,minimum size=4pt] (c) {}; 
			\path (b) ++(120:1) node[circle,fill=black,draw,inner sep=0pt,minimum size=4pt] (d) {}; 
			\node[circle,fill=black,draw,inner sep=0pt,minimum size=4pt] (e) at (1,1) {};
			\node[circle,fill=black,draw,inner sep=0pt,minimum size=4pt] (f) at (1.8,1.4) {};
			\node[circle,fill=black,draw,inner sep=0pt,minimum size=4pt] (g) at (1,1.8) {};
			\node[circle,fill=black,draw,inner sep=0pt,minimum size=4pt] (h) at (1.8,2.2) {};
			\node[circle,fill=black,draw,inner sep=0pt,minimum size=4pt] (i) at (1,2.6) {};
			\node[circle,fill=black,draw,inner sep=0pt,minimum size=4pt] (j) at (1.8,3) {};
			\node[circle,fill=black,draw,inner sep=0pt,minimum size=4pt] (k) at (1,3.4) {};
			\node[circle,fill=black,draw,inner sep=0pt,minimum size=4pt] (l) at (1.8,3.8) {};
			\node[circle,fill=black,draw,inner sep=0pt,minimum size=4pt] (m) at (1,4.2) {};
			\node[circle,fill=black,draw,inner sep=0pt,minimum size=4pt] (n) at (1.8,4.6) {};
			
			\draw[ thick,dotted] (a) -- (b);
			\draw[ thick] (b) to node[midway,inner sep=1pt,outer sep=1pt,minimum size=4pt,fill=white] {$\epsilon$} (c) (b) to node[midway,inner sep=0.5pt,outer sep=0.5pt,minimum size=4pt,fill=white] {$\zeta$} (d) (b) to node[midway,inner sep=1pt,outer sep=1pt,minimum size=4pt,fill=white] {$\eta$} (e) to node[midway,inner sep=1pt,outer sep=1pt,minimum size=4pt,fill=white] {$\alpha$} (f) to node[midway,inner sep=1pt,outer sep=1pt,minimum size=4pt,fill=white] {$\beta$} (g) to node[midway,inner sep=1pt,outer sep=1pt,minimum size=4pt,fill=white] {$\alpha$} (h) to node[midway,inner sep=1pt,outer sep=1pt,minimum size=4pt,fill=white] {$\beta$} (i) to node[midway,inner sep=1pt,outer sep=1pt,minimum size=4pt,fill=white] {$\alpha$} (j) to node[midway,inner sep=1pt,outer sep=1pt,minimum size=4pt,fill=white] {$\beta$} (k) to node[midway,inner sep=1pt,outer sep=1pt,minimum size=4pt,fill=white] {$\alpha$} (l) to node[midway,inner sep=1pt,outer sep=1pt,minimum size=4pt,fill=white] {$\beta$} (m) to node[midway,inner sep=1pt,outer sep=1pt,minimum size=4pt,fill=white] {$\alpha$} (n);
			
			\draw[decoration={brace,amplitude=10pt,mirror},decorate] (2.3,0.9) -- node [midway,below,sloped,yshift=-7pt] {$P_0$} (2.3,3.1);
			
			
			\node[anchor=148.28] at (e) {$z$};
			\node[anchor=north] at (a) {$y$};
			\node[anchor=north] at (b) {$x$};
			\node[anchor=east] at (i) {$x'$};
			\node[anchor=west] at (j) {$y'$};
			\end{tikzpicture}
			\caption{Before shifting.}\label{fig:truncate:initial}
		\end{subfigure}%
		\qquad\qquad\qquad%
		\begin{subfigure}[t]{.2\textwidth}
			\centering
			\begin{tikzpicture}[scale=1.3]
			\node[circle,fill=black,draw,inner sep=0pt,minimum size=4pt] (a) at (0,0) {};
			\node[circle,fill=black,draw,inner sep=0pt,minimum size=4pt] (b) at (1,0) {};
			\path (b) ++(150:1) node[circle,fill=black,draw,inner sep=0pt,minimum size=4pt] (c) {}; 
			\path (b) ++(120:1) node[circle,fill=black,draw,inner sep=0pt,minimum size=4pt] (d) {}; 
			\node[circle,fill=black,draw,inner sep=0pt,minimum size=4pt] (e) at (1,1) {};
			\node[circle,fill=black,draw,inner sep=0pt,minimum size=4pt] (f) at (1.8,1.4) {};
			\node[circle,fill=black,draw,inner sep=0pt,minimum size=4pt] (g) at (1,1.8) {};
			\node[circle,fill=black,draw,inner sep=0pt,minimum size=4pt] (h) at (1.8,2.2) {};
			\node[circle,fill=black,draw,inner sep=0pt,minimum size=4pt] (i) at (1,2.6) {};
			\node[circle,fill=black,draw,inner sep=0pt,minimum size=4pt] (j) at (1.8,3) {};
			\node[circle,fill=black,draw,inner sep=0pt,minimum size=4pt] (k) at (1,3.4) {};
			\node[circle,fill=black,draw,inner sep=0pt,minimum size=4pt] (l) at (1.8,3.8) {};
			\node[circle,fill=black,draw,inner sep=0pt,minimum size=4pt] (m) at (1,4.2) {};
			\node[circle,fill=black,draw,inner sep=0pt,minimum size=4pt] (n) at (1.8,4.6) {};
			
			\draw[ thick] (a) to node[midway,inner sep=1pt,outer sep=1pt,minimum size=4pt,fill=white] {$\epsilon$} (b);
			\draw[dotted, thick] (i) -- (j);
			\draw[ thick] (b) to node[midway,inner sep=0.5pt,outer sep=0.5pt,minimum size=4pt,fill=white] {$\zeta$} (c) (b) to node[midway,inner sep=1pt,outer sep=1pt,minimum size=4pt,fill=white] {$\eta$} (d) (b) to node[midway,inner sep=1pt,outer sep=1pt,minimum size=4pt,fill=white] {$\alpha$} (e) to node[midway,inner sep=1pt,outer sep=1pt,minimum size=4pt,fill=white] {$\beta$} (f) to node[midway,inner sep=1pt,outer sep=1pt,minimum size=4pt,fill=white] {$\alpha$} (g) to node[midway,inner sep=1pt,outer sep=1pt,minimum size=4pt,fill=white] {$\beta$} (h) to node[midway,inner sep=1pt,outer sep=1pt,minimum size=4pt,fill=white] {$\alpha$} (i)  (j) to node[midway,inner sep=1pt,outer sep=1pt,minimum size=4pt,fill=white] {$\beta$} (k) to node[midway,inner sep=1pt,outer sep=1pt,minimum size=4pt,fill=white] {$\alpha$} (l) to node[midway,inner sep=1pt,outer sep=1pt,minimum size=4pt,fill=white] {$\beta$} (m) to node[midway,inner sep=1pt,outer sep=1pt,minimum size=4pt,fill=white] {$\alpha$} (n);
			
			
			
			\node[anchor=148.28] at (e) {$z$};
			\node[anchor=north] at (a) {$y$};
			\node[anchor=north] at (b) {$x$};
			\node[anchor=east] at (i) {$x'$};
			\node[anchor=west] at (j) {$y'$};
			\end{tikzpicture}
			\caption{After shifting.}\label{fig:truncate:shifted}
		\end{subfigure}%
		\caption{Shifting a truncated Vizing chain $F + P_0$.}\label{fig:truncate}
	\end{figure}

	The problem, of course, is that there is no \emph{a priori} upper bound on the length of the alternating path $P$ \ep{except $O(n)$}. In \cite{SV}, Su and Vu tackled this problem as follows. Suppose that $P$ is too long, say, much longer than $\log n$. Pick a random edge $x'y'$ on $P$ within distance $O(\log n)$ from $z$ and truncate the path there. Then let $P_0$ be the truncated path and, instead of shifting the whole chain $F + P$, only shift $F + P_0$. This process is shown in Fig.~\ref{fig:truncate}. In the resulting coloring, the edge $x'y'$ becomes uncolored. If we could use one extra color, say $\Delta + 2$, then we would assign it to $x'y'$. By executing this procedure for all the uncolored edges, one eventually obtains a $(\Delta+2)$-edge-coloring of the entire graph that is proper except that there might be adjacent edges of color $\Delta + 2$. However, Su and Vu showed in \cite{SV} that with high probability, no two edges of color $\Delta + 2$ will be adjacent, and thus one actually gets a proper $(\Delta+2)$-edge-coloring.
	
\begin{figure}[b]
		\centering
		\begin{tikzpicture}[scale=1.3]
		\node[circle,fill=black,draw,inner sep=0pt,minimum size=4pt] (a) at (0,0) {};
		\node[circle,fill=black,draw,inner sep=0pt,minimum size=4pt] (b) at (1,0) {};
		\path (b) ++(150:1) node[circle,fill=black,draw,inner sep=0pt,minimum size=4pt] (c) {};
		\path (b) ++(120:1) node[circle,fill=black,draw,inner sep=0pt,minimum size=4pt] (d) {}; 
		\node[circle,fill=black,draw,inner sep=0pt,minimum size=4pt] (e) at (1,1) {};
		\node[circle,fill=black,draw,inner sep=0pt,minimum size=4pt] (f) at (1.8,1.4) {};
		\node[circle,fill=black,draw,inner sep=0pt,minimum size=4pt] (g) at (1,1.8) {};
		\node[circle,fill=black,draw,inner sep=0pt,minimum size=4pt] (h) at (1.8,2.2) {};
		\node[circle,fill=black,draw,inner sep=0pt,minimum size=4pt] (i) at (1,2.6) {};
		\node[circle,fill=black,draw,inner sep=0pt,minimum size=4pt] (l) at (1,3.6) {};
		
		
		\path (i) ++(60:1) node[circle,fill=black,draw,inner sep=0pt,minimum size=4pt] (m) {};
		\path (i) ++(30:1) node[circle,fill=black,draw,inner sep=0pt,minimum size=4pt] (n) {};
		\path (i) ++(0:1) node[circle,fill=black,draw,inner sep=0pt,minimum size=4pt] (o) {};
		
		\node[circle,fill=black,draw,inner sep=0pt,minimum size=4pt] (p) at (2.4,3.4) {};
		\node[circle,fill=black,draw,inner sep=0pt,minimum size=4pt] (q) at (2.8,2.6) {};
		\node[circle,fill=black,draw,inner sep=0pt,minimum size=4pt] (r) at (3.2,3.4) {};
		\node[circle,fill=black,draw,inner sep=0pt,minimum size=4pt] (s) at (3.6,2.6) {};
		\node[circle,fill=black,draw,inner sep=0pt,minimum size=4pt] (t) at (4,3.4) {};
		\node[circle,fill=black,draw,inner sep=0pt,minimum size=4pt] (u) at (4.4,2.6) {};
		
		
		
		\draw[ thick,dotted] (a) -- (b);
		\draw[ thick] (b) -- (c) (b) -- (d) (b) -- (e) to node[midway,inner sep=1pt,outer sep=1pt,minimum size=4pt,fill=white] {$\alpha$} (f) to node[midway,inner sep=1pt,outer sep=1pt,minimum size=4pt,fill=white] {$\beta$} (g) to node[midway,inner sep=1pt,outer sep=1pt,minimum size=4pt,fill=white] {$\alpha$} (h) to node[midway,inner sep=1pt,outer sep=1pt,minimum size=4pt,fill=white] {$\beta$} (i) to node[midway,inner sep=1pt,outer sep=1pt,minimum size=4pt,fill=white] {$\alpha$} (l);
		
		\draw[ thick] (i) -- (m) (i) -- (n) (i) -- (o) to node[midway,inner sep=1pt,outer sep=1pt,minimum size=4pt,fill=white] {$\gamma$} (p) to node[midway,inner sep=1pt,outer sep=1pt,minimum size=4pt,fill=white] {$\delta$} (q) to node[midway,inner sep=1pt,outer sep=1pt,minimum size=4pt,fill=white] {$\gamma$} (r) to node[midway,inner sep=1pt,outer sep=1pt,minimum size=4pt,fill=white] {$\delta$} (s) to node[midway,inner sep=1pt,outer sep=1pt,minimum size=4pt,fill=white] {$\gamma$} (t) to node[midway,inner sep=1pt,outer sep=1pt,minimum size=4pt,fill=white] {$\delta$} (u);
		
		\node[anchor=148.28] at (e) {$z$};
		\node[anchor=north] at (a) {$y$};
		\node[anchor=north] at (b) {$x$};
		\node[anchor=east] at (i) {$x'$};
		\node[anchor=east] at (l) {$y'$};
		\node[anchor=166.71] at (o) {$z'$};
		
		
		\draw[decoration={brace,amplitude=10pt},decorate] (0.6,1) -- node [midway,above,sloped,yshift=7pt] {$P_0$} (0.6,3.7);
		
		\draw[decoration={brace,amplitude=10pt},decorate] (2,3.4) -- node [midway,above,sloped,yshift=7pt] {$P'$} (4.5,3.4);
		
		\end{tikzpicture}
		\caption{A two-step Vizing chain $F + P_0 + F' + P'$.}\label{fig:twostep}
	\end{figure}
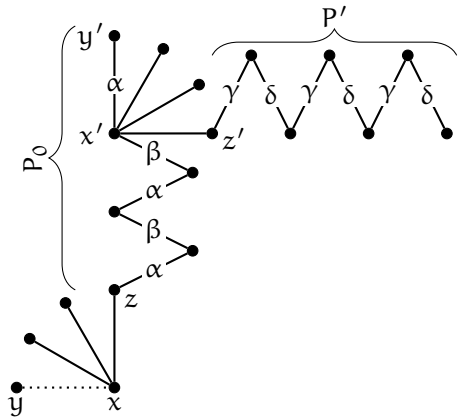

	In \cite{GP}, Greb\'\i{}k and Pikhurko encountered a similar problem, but in the context of $(\Delta + 1)$-edge-coloring \ep{as we already mentioned in \S\ref{subsec:intro_comb}, Greb\'\i{}k and Pikhurko's goal was to prove a version of Vizing's theorem for infinite graphs on probability spaces}. Their solution also involved picking a random edge $x'y'$ on $P$, truncating the path there, and then shifting the chain $F + P_0$. But since they could not use an extra color for the edge $x'y'$, they instead built a \emph{new} Vizing chain $H' = F' + P'$ starting at the edge $x'y'$, and shifted the whole combination of $F + P_0 + F' + P'$ \ep{see Fig.~\ref{fig:twostep}}. Let us call the sequence $F + P_0 + F' + P'$ a \emph{two-step Vizing chain}. This two-step approach can be used to build an augmenting subgraph with only $O(\sqrt{n})$ edges, and while our ultimate goal is to reduce $O(\sqrt{n})$ further to $O((\log n)^2)$, it will perhaps be instructive if we explain first how to get the bound $O(\sqrt{n})$.
	
	If the length of $P$ is $O(\sqrt{n})$, then we are done, so assume that $P$ is much longer than $\sqrt{n}$. The vertex $x'$ is chosen at random from among the first $\Theta(\sqrt{n})$ vertices of $P$, and we want to show that the expected length of $P'$, relative to this random choice of $x'$, is $O(\sqrt{n})$. To this end, we fix an arbitrary vertex $u \in V(G)$ and ask, what is the probability that $P'$ passes through $u$? Suppose that the path $P'$ starts at a vertex $z'$ that is adjacent to $x'$ and is $\gamma\delta$-alternating \ep{as is shown in Fig.~\ref{fig:twostep}}. If we start at $u$ and follow the edges colored $\gamma$ and $\delta$, then we would stop precisely when we hit $z'$. Thus, given $u$, we can uniquely determine $x'$ by specifying the colors $\gamma$ and $\delta$ \ep{$(\Delta+1)^2$ options} and the location of $x'$ among the neighbors of $z'$ \ep{$\Delta$ options}. In other words, given $u$, there are at most $O(\Delta^3) = O(1)$ possibilities for $x'$. Since $x'$ is chosen randomly from a set of $\Theta(\sqrt{n})$ candidates, we conclude that
	\[
	\P\left[\text{$P'$ passes through $u$}\right] \,\leq\, O(1/\sqrt{n}).
	\]
	Summing over all $u \in V(G)$ yields $\mathbb{E}[|P'|] \leq O(n \cdot 1/\sqrt{n}) = O(\sqrt{n})$, as desired.
	
	The argument in the previous paragraph has a flaw that, even in this informal overview, we cannot overlook. Namely, the path $P'$ may intersect the fan $F$ or the path $P_0$. Because of that, some of the edges of $P'$ may be colored differently before and after the chain $F + P_0$ is shifted. This means that, in principle, even if we know $u$, $\delta$, and $\gamma$, we cannot locate $z'$ by tracing the $\gamma\delta$-alternating path starting at $u$ unless we also know what the coloring looks like after shifting $F + P_0$. To circumvent this issue, Greb\'\i{}k and Pikhurko showed that, roughly speaking, a careful choice of $P'$ ensures that $P'$ is disjoint from $F + P_0$ with high probability. It is unlikely that $P'$ intersects $F$---essentially because $F$ is small, namely $|F| \leq \Delta = O(1)$. The real challenge is to make $P'$ \ep{edge-}disjoint from $P_0$. Since every edge in $P_0$ is colored $\alpha$ or $\beta$, it would be ideal if we could simply pick $P'$ so that $\set{\gamma, \delta} \cap \set{\alpha, \beta} = \0$. Unfortunately, this is not always possible. Thankfully, Greb\'\i{}k and Pikhurko showed that $P'$ can always be chosen so that either $\set{\gamma, \delta} \cap \set{\alpha, \beta} = \0$ or else, $\set{\gamma, \delta} = \set{\alpha, \beta}$ \ep{see Lemma~\ref{lemma:fancy_fan} below}. In the former case, $P'$ and $P_0$ are disjoint, so we are done. In the latter case, $P'$ and $P_0$ cannot intersect transversally: if $P' \cap P_0 \neq \0$, then $P' \supseteq P_0$, and one can show that this is unlikely.
	
\begin{figure}[h]
		\centering
		\begin{tikzpicture}[scale=1.3]
		\node[circle,fill=black,draw,inner sep=0pt,minimum size=4pt] (a) at (0,0) {};
		\node[circle,fill=black,draw,inner sep=0pt,minimum size=4pt] (b) at (1,0) {};
		\path (b) ++(150:1) node[circle,fill=black,draw,inner sep=0pt,minimum size=4pt] (c) {};
		\path (b) ++(120:1) node[circle,fill=black,draw,inner sep=0pt,minimum size=4pt] (d) {};
		\node[circle,fill=black,draw,inner sep=0pt,minimum size=4pt] (e) at (1,1) {};
		\node[circle,fill=black,draw,inner sep=0pt,minimum size=4pt] (f) at (1.8,1.4) {};
		\node[circle,fill=black,draw,inner sep=0pt,minimum size=4pt] (g) at (1,1.8) {};
		\node[circle,fill=black,draw,inner sep=0pt,minimum size=4pt] (h) at (1.8,2.2) {};
		\node[circle,fill=black,draw,inner sep=0pt,minimum size=4pt] (i) at (1,2.6) {};
		\node[circle,fill=black,draw,inner sep=0pt,minimum size=4pt] (j) at (1.8,3) {};
		\node[circle,fill=black,draw,inner sep=0pt,minimum size=4pt] (k) at (1,3.4) {};
		\node[circle,fill=black,draw,inner sep=0pt,minimum size=4pt] (l) at (1,4.4) {};
		
		\path (k) ++(60:1) node[circle,fill=black,draw,inner sep=0pt,minimum size=4pt] (m) {};
		\path (k) ++(30:1) node[circle,fill=black,draw,inner sep=0pt,minimum size=4pt] (n) {};
		\path (k) ++(0:1) node[circle,fill=black,draw,inner sep=0pt,minimum size=4pt] (o) {};
		
		\node[circle,fill=black,draw,inner sep=0pt,minimum size=4pt] (p) at (2.4,4.2) {};
		\node[circle,fill=black,draw,inner sep=0pt,minimum size=4pt] (q) at (2.8,3.4) {};
		\node[circle,fill=black,draw,inner sep=0pt,minimum size=4pt] (r) at (3.2,4.2) {};
		\node[circle,fill=black,draw,inner sep=0pt,minimum size=4pt] (s) at (3.6,3.4) {};
		\node[circle,fill=black,draw,inner sep=0pt,minimum size=4pt] (t) at (4,4.2) {};

		\node[circle,fill=black,draw,inner sep=0pt,minimum size=4pt] (v) at (5,4.2) {};
		\path (t) ++(-30:1) node[circle,fill=black,draw,inner sep=0pt,minimum size=4pt] (w) {};
		\path (t) ++(-60:1) node[circle,fill=black,draw,inner sep=0pt,minimum size=4pt] (x) {};
		\path (t) ++(-90:1) node[circle,fill=black,draw,inner sep=0pt,minimum size=4pt] (y) {};

		\node[circle,fill=black,draw,inner sep=0pt,minimum size=4pt] (z) at (4.8,2.8) {};
		\node[circle,fill=black,draw,inner sep=0pt,minimum size=4pt] (aa) at (4,2.4) {};
		\node[circle,fill=black,draw,inner sep=0pt,minimum size=4pt] (ab) at (4.8,2) {};
		\node[circle,fill=black,draw,inner sep=0pt,minimum size=4pt] (ac) at (4,1.6) {};
		\node[circle,fill=black,draw,inner sep=0pt,minimum size=4pt] (ad) at (4.8,1.2) {};
		\node[circle,fill=black,draw,inner sep=0pt,minimum size=4pt] (ae) at (4.8,0.2) {};
		
		\path (ad) ++(-120:1) node[circle,fill=black,draw,inner sep=0pt,minimum size=4pt] (ag) {};
		\path (ad) ++(-150:1) node[circle,fill=black,draw,inner sep=0pt,minimum size=4pt] (ah) {};
		\path (ad) ++(-180:1) node[circle,fill=black,draw,inner sep=0pt,minimum size=4pt] (ai) {};
		
		\node[circle,fill=black,draw,inner sep=0pt,minimum size=4pt] (aj) at (3.4,0.4) {};
		\node[circle,fill=black,draw,inner sep=0pt,minimum size=4pt] (ak) at (3,1.2) {};
		\node[circle,fill=black,draw,inner sep=0pt,minimum size=4pt] (al) at (2.6,0.4) {};
		\node[circle,fill=black,draw,inner sep=0pt,minimum size=4pt] (am) at (2.2,1.2) {};

		\draw[ thick,dotted] (a) -- (b);
		\draw[ thick] (b) -- (c) (b) -- (d) (b) -- (e) to node[midway,inner sep=1pt,outer sep=1pt,minimum size=4pt,fill=white] {$\alpha_0$} (f) to node[midway,inner sep=1pt,outer sep=1pt,minimum size=4pt,fill=white] {$\beta_0$} (g) to node[midway,inner sep=1pt,outer sep=1pt,minimum size=4pt,fill=white] {$\alpha_0$} (h) to node[midway,inner sep=1pt,outer sep=1pt,minimum size=4pt,fill=white] {$\beta_0$} (i) to node[midway,inner sep=1pt,outer sep=1pt,minimum size=4pt,fill=white] {$\alpha_0$} (j) to node[midway,inner sep=1pt,outer sep=1pt,minimum size=4pt,fill=white] {$\beta_0$} (k) to node[midway,inner sep=1pt,outer sep=1pt,minimum size=4pt,fill=white] {$\alpha_0$} (l);
		
		\draw[ thick] (k) -- (m) (k) -- (n) (k) -- (o) to node[midway,inner sep=1pt,outer sep=1pt,minimum size=4pt,fill=white] {$\alpha_1$} (p) to node[midway,inner sep=1pt,outer sep=1pt,minimum size=4pt,fill=white] {$\beta_1$} (q) to node[midway,inner sep=1pt,outer sep=1pt,minimum size=4pt,fill=white] {$\alpha_1$} (r) to node[midway,inner sep=1pt,outer sep=1pt,minimum size=4pt,fill=white] {$\beta_1$} (s) to node[midway,inner sep=1pt,outer sep=1pt,minimum size=4pt,fill=white] {$\alpha_1$} (t) to node[midway,inner sep=1pt,outer sep=1pt,minimum size=4pt,fill=white] {$\beta_1$} (v);
		
		\draw[ thick] (t) -- (w) (t) -- (x) (t) -- (y) to node[midway,inner sep=1pt,outer sep=1pt,minimum size=4pt,fill=white] {$\alpha_2$} (z) to node[midway,inner sep=1pt,outer sep=1pt,minimum size=4pt,fill=white] {$\beta_2$} (aa) to node[midway,inner sep=1pt,outer sep=1pt,minimum size=4pt,fill=white] {$\alpha_2$} (ab) to node[midway,inner sep=1pt,outer sep=1pt,minimum size=4pt,fill=white] {$\beta_2$} (ac) to node[midway,inner sep=1pt,outer sep=1pt,minimum size=4pt,fill=white] {$\alpha_2$} (ad) to node[midway,inner sep=1pt,outer sep=1pt,minimum size=4pt,fill=white] {$\beta_2$} (ae);
		
		\draw[thick] (ad) -- (ag) (ad) -- (ah) (ad) -- (ai) to node[midway,inner sep=1pt,outer sep=1pt,minimum size=4pt,fill=white] {$\alpha_3$} (aj) to node[midway,inner sep=1pt,outer sep=1pt,minimum size=4pt,fill=white] {$\beta_3$} (ak) to node[midway,inner sep=1pt,outer sep=1pt,minimum size=4pt,fill=white] {$\alpha_3$} (al) to node[midway,inner sep=1pt,outer sep=1pt,minimum size=4pt,fill=white] {$\beta_3$} (am);
		
		\node[anchor=148.28] at (e) {$z_0$};
		\node[anchor=north] at (a) {$y_0$};
		\node[anchor=north] at (b) {$x_0$};
		
		\node[anchor=121.71] at (o) {$z_1$};
		\node[anchor=east] at (l) {$y_1$};
		\node[anchor=east] at (k) {$x_1$};
		
		\node[anchor=31.71] at (y) {$z_2$};
		\node[anchor=south] at (v) {$y_2$};
		\node[anchor=south] at (t) {$x_2$};
		
		\node[anchor=-58.28] at (ai) {$z_3$};
		\node[anchor=west] at (ae) {$y_3$};
		\node[anchor=west] at (ad) {$x_3$};
		\end{tikzpicture}
		\caption{A multi-step Vizing chain.}\label{fig:multistep}
	\end{figure}
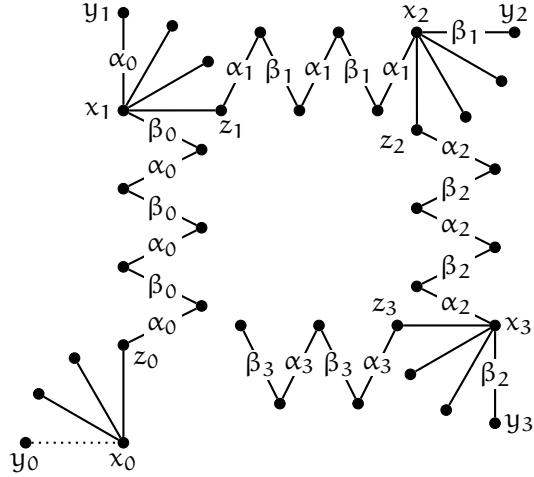

	In this paper we take a natural next step and consider \emph{multi-step Vizing chains} \ep{see Fig.~\ref{fig:multistep}}. Fix parameters $\ell = a \log n$ and $T = b \log n$, where $a$ and $b$ are positive constants with $a \gg b$. Let $x_0 \defeq x$ and $y_0 \defeq y$ and build a Vizing chain $F_0 + Q_0$ for $x_0y_0$. Here $Q_0$ is an $\alpha_0\beta_0$-alternating path that starts at a vertex $z_0$ adjacent to $x_0$. If $|Q_0| \leq \ell$, then stop. Otherwise, pick a random edge $x_1y_1$ on $Q_0$ at distance at most $\ell$ from $z_0$, truncate $Q_0$ at that edge, and let $P_0$ be the truncated path. After shifting $F_0 + P_0$, the edge $x_1y_1$ becomes uncolored, and we repeat the process with $x_1y_1$ in place of $x_0y_0$, obtaining a fan $F_1$ and an $\alpha_1\beta_1$-alternating path $Q_1$ starting at a vertex $z_1$ adjacent to $x_1$. Again, if $|Q_1| \leq \ell$, then we stop, and otherwise we pick a random edge $x_2y_2$ on $Q_1$, truncate $Q_1$ at that edge, and let $P_1$ be the truncated path. Then we shift $F_1 + P_1$, and repeat the process with $x_2y_2$ in place of $x_1y_1$, etc. If this process terminates within $T$ steps, then it produces an augmenting subgraph with at most $O(T\ell) = O((\log n)^2)$ edges, as desired. Hence, our goal becomes to show that the process terminates within $T$ steps with positive probability.
	
	The basic calculation here is similar to the one in the two-step case. Namely, we fix a vertex $u \in V(G)$ and bound the probability that the process runs for at least $T$ steps and the path $P_T$ passes through $u$. If we start at $u$ and follow the edges colored $\alpha_T$ and $\beta_T$, we will stop at $z_T$. Since $x_T$ is adjacent to $z_T$, we can determine $x_T$ by specifying the colors $\alpha_T$ and $\beta_T$ and the location of $x_T$ among the neighbors of $z_T$. Then we can trace the $\alpha_{T-1}\beta_{T-1}$-alternating path starting at $x_T$ to find $z_{T-1}$, and then locate $x_{T-1}$ among the neighbors of $z_{T-1}$. Continuing in this manner, we see that the total number of sequences $(x_1, \ldots, x_T)$ for which $P_T$ passes through $u$ is at most $\Delta^{O(T)}$. Since each $x_i$ is chosen among $\ell$ candidates, we conclude that, for some positive constant $C$,
	\[
	\P\left[\text{$P_T$ passes through $u$}\right] \,\leq\, \frac{\Delta^{O(T)}}{\ell^T} \,\leq\, \exp(-C \log n \log \log n).
	\]
	If the process runs for at least $T$ steps, then $P_T$ has to pass through \emph{some} vertex, so
	\[
	\P\left[\text{the process runs for at least $T$ steps}\right] \,\leq\, n \cdot \exp(-C \log n \log \log n) \,=\, o(1).
	\]
	
	Again, the above reasoning only really works if each $P_i$ is disjoint from $F_0 + P_0 + \cdots + F_{i-1} + P_{i-1}$. As in the two-step case, making $P_i$ disjoint from the fans $F_0$, \ldots, $F_{i-1}$ is a bit easier, so let us focus on the problem of making $P_i$ disjoint from $P_0$, \ldots, $P_{i-1}$. For $i > j$, let $\mathcal{E}_{ij}$ be the event that the first intersection that occurs during the process is between $P_i$ and $P_j$. We need to bound the probability of $\mathcal{E}_{ij}$. By using the same methods as Greb\'\i{}k and Pikhurko, i.e., by making sure that either $\set{\alpha_i, \beta_i} \cap \set{\alpha_{i-1}, \beta_{i-1}} = \0$ or $\set{\alpha_i, \beta_i} = \set{\alpha_{i-1}, \beta_{i-1}}$, we can eliminate the possibility that $j = i - 1$, so let us assume that $j \leq i - 2$. Fix the outcome of the procedure up to step $j$, so that we have already settled on the vertices $x_0$, $x_1$, \ldots, $x_j$ but not on $x_{j+1}$. At this point, we already know the path $Q_j$, but we have not yet truncated it. Pick any vertex $u$ on $Q_j$ within distance $\ell$ from $z_j$ and let $\mathcal{E}_{ij,u}$ be the event that $\mathcal{E}_{ij}$ happens and the first intersection point between $P_i$ and $P_j$ is $u$. The key observation is that if $\mathcal{E}_{ij,u}$ happens, then, since $u$ is the \emph{first} intersection point, the calculations in the preceding paragraph \emph{can} be used for $u$. That is, starting at $u$, we can trace the $\alpha_i\beta_i$-alternating path back to $z_i$, choose $x_i$ from the neighbors of $z_i$, trace the $\alpha_{i-1}\beta_{i-1}$-alternating path to $z_{i-1}$, choose $x_{i-1}$ from the neighbors of $z_{i-1}$, etc. In this way, we see that there is a positive constant $c$ such that  the number of sequences $(x_{j+1}, \ldots, x_i)$ for which $\mathcal{E}_{ij,u}$ holds is at most $\Delta^{c (i-j)}$. Therefore,
	\[
	\P\left[\mathcal{E}_{ij, u}\right] \,\leq\, \frac{\Delta^{c(i-j)}}{\ell^{i - j}}.
	\]
	There are $\ell$ options for $u$, so, by the union bound,
	\[
	\P\left[\mathcal{E}_{ij}\right] \,\leq\, \ell \cdot \frac{\Delta^{c(i-j)}}{\ell^{i - j}} \,=\, \frac{\Delta^{c(i-j)}}{\ell^{i - j -1}}.
	\]	
	Applying the union bound a second time and using the assumption $j \leq i - 2$, we get
	\[
	\P\left[\text{the paths $P_0$, \ldots, $P_T$ are not disjoint}\right] \,\leq\, \sum_{j = 0}^{T-2} \sum_{i = j+2}^T \frac{\Delta^{c(i-j)}}{\ell^{i - j -1}} \,\leq\, \Delta^{c}\sum_{j = 0}^{T-2}\sum_{k = 1}^\infty \frac{\Delta^{ck}}{\ell^{k}} \,=\, \frac{\Delta^{2c} (T-1)}{\ell - \Delta^c}.
	\]
	Remembering that $\ell = a \log n$ and $T = b \log n$, we can make the last expression as small as desired by increasing the ratio $a/b$. Notice that $O((\log n)^2)$ is a natural barrier for this argument, since both $\ell$ and $T$ must be of the order $\Theta(\log n)$. \ep{With more careful computations, it is possible to shave off a $\mathrm{poly}(\log \log n)$-factor from this bound.} It remains an interesting open problem to see if one can actually reduce the bound to $O(\log n)$ using a different method.
	
	The rest of the paper is organized as follows. In \S\ref{sec:prelim}, we introduce the terminology needed for the proof, such as the notion of shifting a chain. We also prove Lemmas \ref{lemma:simple_fan} and \ref{lemma:fancy_fan}, the first of which essentially asserts the existence of a Vizing chain, while the second one comes from the work of Greb\'\i{}k and Pikhurko and allows one to find $\gamma\delta$-alternating paths with $\set{\gamma, \delta} \cap \set{\alpha, \beta} = \0$ or $\set{\gamma, \delta} = \set{\alpha, \beta}$. Then, in \S\ref{sec:const}, we describe our Multi-Step Vizing Algorithm in detail and verify some of its properties. After that, in \S\ref{sec:analysis}, we undertake the probabilistic analysis of the Multi-Step Vizing Algorithm and finish the proof of Theorem~\ref{theo:ext}.

	\section{Chains, paths, and fans}\label{sec:prelim}
	
	
	\subsection{Basics}
	
	For the remainder of the paper we fix an $n$-vertex graph $G$ of maximum degree $\Delta$. For brevity, set $V \defeq V(G)$ and $E \defeq E(G)$. As usual, we say that two edges of $G$ are \emphd{adjacent}, or \emphd{neighbors} of each other, if they are distinct and have a common endpoint. Throughout, we shall use the phrase ``\ep{partial} coloring'' to mean a \ep{partial} $(\Delta+1)$-edge-coloring of $G$. For convenience, if $\phi$ is a partial coloring and $e \in E \setminus \dom(\phi)$ is an uncolored edge, then we write $\phi(e) = \blank$, where $\blank$ is a special symbol distinct from every color.
	
	\subsection{Missing colors and happy edges}
	
	Given a proper partial coloring $\phi$ and a vertex $x \in V$, we let $M(\phi, x) \subseteq [\Delta + 1]$ be the set of all colors \emphd{missing} at $x$ in the coloring $\phi$, that is,
	$M(\phi, x) \defeq [\Delta + 1] \setminus \set{\phi(xy) \,:\, xy \in E}$.
	Note that $|M(\phi, x)|$ is equal to $\Delta + 1 - \deg_G(x)$ plus the number of uncolored edges incident to $x$; in particular, the set $M(\phi, x)$ is always nonempty. We say that an uncolored edge $e = xy$ is \emphd{$\phi$-happy} if $M(\phi, x) \cap M(\phi, y) \neq \0$. If $e = xy$ is $\phi$-happy and $\alpha$ is any color in $M(\phi, x) \cap M(\phi, y)$, then we can extend $\phi$ by assigning to $e$ the color $\alpha$. In other words, in this case the subgraph of $G$ with vertex set $\set{x,y}$ and edge set $\set{e}$ is augmenting for $\phi$. 
	
	\subsection{Shifts and chains}
	
	Given a proper partial coloring, we wish to ``locally'' modify it in order to create a happy edge \ep{which then would allow us to enlarge the domain of the coloring}. In order to achieve this, we will be ``shifting'' colors from colored edges to their uncolored neighbors. Let $e_0$, $e_1 \in E$ be two adjacent edges and let $\phi$ be a proper partial coloring. Define a 
	coloring $\Shift(\phi, e_0, e_1)$ by setting, for all $e \in E$,
	\[
	\Shift(\phi, e_0, e_1)(e) \,\defeq\, \begin{cases}
	\phi(e_1) &\text{if } e = e_0;\\
	\blank &\text{if } e = e_1;\\
	\phi(e) &\text{otherwise}.
	\end{cases}
	\]
	In other words, $\Shift(\phi, e_0, e_1)$ ``shifts'' the color from $e_1$ to $e_0$, leaves $e_1$ uncolored, and keeps the coloring of the rest of the edges unchanged. The pair $(e_0, e_1)$ is said to be \emphd{$\phi$-shiftable} if
	\begin{itemize}
		\item $e_0 \not \in \dom(\phi)$, $e_1 \in \dom(\phi)$, and
		\item the coloring $\Shift(\phi, e_0, e_1)$ is proper.
	\end{itemize}
	The latter condition is equivalent to saying that $\phi(e_1) \in M(\phi, x)$, where $x$ is the unique vertex that belongs to $e_0$ but not to $e_1$. Note that if a pair $(e_0, e_1)$ is $\phi$-shiftable, then the pair $(e_1, e_0)$ is $\Shift(\phi, e_0, e_1)$-shiftable and 
	\begin{equation}\label{eq:reverse}
	\Shift(\Shift(\phi, e_0, e_1), e_1, e_0) \,=\, \phi.
	\end{equation}
	
\begin{figure}[t]
		\centering
		\begin{tikzpicture}
		\begin{scope}
		\node[circle,fill=black,draw,inner sep=0pt,minimum size=4pt] (a) at (0,0) {};
		\node[circle,fill=black,draw,inner sep=0pt,minimum size=4pt] (b) at (-1.25,0.625) {};
		\node[circle,fill=black,draw,inner sep=0pt,minimum size=4pt] (c) at (0,1.625) {};
		\node[circle,fill=black,draw,inner sep=0pt,minimum size=4pt] (d) at (1.25,0.625) {};
		\node[circle,fill=black,draw,inner sep=0pt,minimum size=4pt] (e) at (2.5,0) {};
		\node[circle,fill=black,draw,inner sep=0pt,minimum size=4pt] (f) at (3.75,0.625) {};
		\node[circle,fill=black,draw,inner sep=0pt,minimum size=4pt] (g) at (3.75,-1) {};
		\node[circle,fill=black,draw,inner sep=0pt,minimum size=4pt] (h) at (5,0) {};
		
		\draw[thick,dotted] (a) to node[midway,inner sep=0pt,minimum size=4pt] (i) {} (b);
		\draw[thick] (a) to node[midway,inner sep=1pt,outer sep=1pt,minimum size=4pt,fill=white] (j) {$\alpha$} (c);
		\draw[thick] (a) to node[midway,inner sep=1pt,outer sep=1pt,minimum size=4pt,fill=white] (k) {$\beta$} (d);
		\draw[thick] (d) to node[midway,inner sep=1pt,outer sep=1pt,minimum size=4pt,fill=white] (l) {$\gamma$} (e);
		\draw[thick] (f) to node[midway,inner sep=1pt,outer sep=1pt,minimum size=4pt,fill=white] (m) {$\delta$} (e);
		\draw[thick] (f) to node[midway,inner sep=1pt,outer sep=1pt,minimum size=4pt,fill=white] (n) {$\epsilon$} (g);
		\draw[thick] (f) to node[midway,inner sep=1pt,outer sep=1pt,minimum size=4pt,fill=white] (o) {$\zeta$} (h);
		
		\draw[-{Stealth[length=1.6mm]}] (j) to[bend right] (i);
		\draw[-{Stealth[length=1.6mm]}] (k) to[bend right] (j);
		\draw[-{Stealth[length=1.6mm]}] (l) to[bend left] (k);
		\draw[-{Stealth[length=1.6mm]}] (m) to[bend right] (l);
		\draw[-{Stealth[length=1.6mm]}] (n) to[bend left] (m);
		\draw[-{Stealth[length=1.6mm]}] (o) to[bend left] (n);
		\end{scope}
		
		\draw[-{Stealth[length=1.6mm]},decoration = {snake,pre length=3pt,post length=7pt,},decorate] (5.4,0.3125) -- (6.35,0.3125);
		
		\begin{scope}[xshift=8cm]
		\node[circle,fill=black,draw,inner sep=0pt,minimum size=4pt] (a) at (0,0) {};
		\node[circle,fill=black,draw,inner sep=0pt,minimum size=4pt] (b) at (-1.25,0.625) {};
		\node[circle,fill=black,draw,inner sep=0pt,minimum size=4pt] (c) at (0,1.625) {};
		\node[circle,fill=black,draw,inner sep=0pt,minimum size=4pt] (d) at (1.25,0.625) {};
		\node[circle,fill=black,draw,inner sep=0pt,minimum size=4pt] (e) at (2.5,0) {};
		\node[circle,fill=black,draw,inner sep=0pt,minimum size=4pt] (f) at (3.75,0.625) {};
		\node[circle,fill=black,draw,inner sep=0pt,minimum size=4pt] (g) at (3.75,-1) {};
		\node[circle,fill=black,draw,inner sep=0pt,minimum size=4pt] (h) at (5,0) {};
		
		\draw[thick] (a) to node[midway,inner sep=1pt,outer sep=1pt,minimum size=4pt,fill=white] (i) {$\alpha$} (b);
		\draw[thick] (a) to node[midway,inner sep=1pt,outer sep=1pt,minimum size=4pt,fill=white] (j) {$\beta$} (c);
		\draw[thick] (a) to node[midway,inner sep=1pt,outer sep=1pt,minimum size=4pt,fill=white] (k) {$\gamma$} (d);
		\draw[thick] (d) to node[midway,inner sep=1pt,outer sep=1pt,minimum size=4pt,fill=white] (l) {$\delta$} (e);
		\draw[thick] (f) to node[midway,inner sep=1pt,outer sep=1pt,minimum size=4pt,fill=white] (m) {$\epsilon$} (e);
		\draw[thick] (f) to node[midway,inner sep=1pt,outer sep=1pt,minimum size=4pt,fill=white] (n) {$\zeta$} (g);
		\draw[thick,dotted] (f) to node[midway,inner sep=0pt,minimum size=4pt] (o) {} (h);
		\end{scope}
		\end{tikzpicture}
		\caption{Shifting a coloring along a chain.}\label{fig:shift}
	\end{figure}
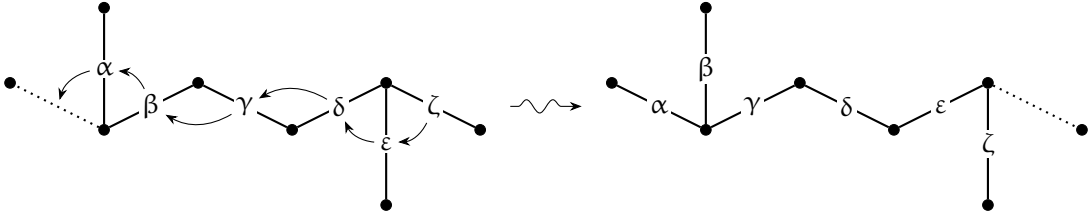

	A \emphd{chain} of length $\ell \geq 1$ is a sequence $C = (e_0, \ldots, e_{\ell-1})$ of edges such that $e_i$ and $e_{i+1}$ are adjacent, for all $0 \leq i \leq \ell-2$. Note that we do not require the edges in a chain to be pairwise distinct; however, any two \emph{consecutive} edges must be distinct since they are adjacent. We define the \emphd{edge set} $E(C)$ and the \emphd{vertex set} $V(C)$ of a chain $C = (e_0, \ldots, e_{\ell-1})$ by setting $E(C) \defeq \set{e_0, \ldots, e_{\ell-1}}$ and making $V(C)$ be the set of all the endpoints of the edges in $C$. We also let $H(C)$ be the subgraph of $G$ with vertex set $V(C)$ and edge set $E(C)$. It is clear that the graph $H(C)$ is connected. Let $\start(C) \defeq e_0$ and $\eend(C) \defeq e_{\ell-1}$ denote the first and the last edges of $C$, respectively, and let $\length(C) \defeq \ell$ be the length of $C$. Given a partial coloring $\phi$, we define $\Shift(\phi, C)$ by iteratively shifting the colors from $e_1$ to $e_0$, then from $e_2$ to $e_1$, etc.~\ep{see Fig.~\ref{fig:shift}}. More precisely, we consider the following recursively defined sequence of colorings:
	\begin{align*}
	\Shift_0(\phi, C) \,&\defeq\, \phi;\\
	\Shift_{i+1}(\phi, C) \,&\defeq\, \Shift(\Shift_{i}(\phi, C), e_i, e_{i+1}) \quad \text{for all } 0 \leq i \leq \ell-2,
	\end{align*}
	and set $\Shift(\phi, C) \defeq \Shift_{\ell - 1}(\phi, C)$. We say that $C$ is \emphd{$\phi$-shiftable} if for all $0 \leq i \leq \ell - 2$, the pair $(e_i, e_{i+1})$ is $\Shift_{i}(\phi, C)$-shiftable. Given a chain $C = (e_0, \ldots, e_{\ell - 1})$, we let $C^\ast \defeq (e_{\ell - 1}, \ldots, e_0)$ denote the \emphd{reverse} of $C$. Repeated applications of \eqref{eq:reverse} show that if $C$ is $\phi$-shiftable, then $C^\ast$ is $\Shift(\phi, C)$-shiftable and $\Shift(\Shift(\phi, C), C^\ast) = \phi$.
	
	\begin{fact}\label{fact:shift}
		Let $\phi$ be a proper partial coloring and let $C$ be a chain. If $C$ is $\phi$-shiftable, then the following statements are valid:
		\begin{enumerate}[label=\ep{\normalfont\roman*}]
			\item\label{item:start} the unique edge in $E(C)$ that is uncolored by $\phi$ is $\start(C)$;
			
			\item\label{item:proper} $\Shift(\phi, C)$ is a proper partial coloring that agrees with $\phi$ on $E \setminus E(C)$;
			
			\item\label{item:end} the unique edge in $E(C)$ that is uncolored by $\Shift(\phi, C)$ is $\eend(C)$.
		\end{enumerate}
	\end{fact}
	\begin{scproof}
		Let $C = (e_0, \ldots, e_{\ell - 1})$, where $\ell= \length(C)$. A straightforward induction shows that for each $0 \leq i \leq \ell - 1$, the coloring $\Shift_i(\phi, C)$ is proper, which, in particular, yields \ref{item:proper}. To see \ref{item:start}, note first that $e_0$ must be uncolored by $\phi$ since the pair $(e_0, e_1)$ is $\phi$-shiftable. Suppose now that $e \in E(C)$ is a $\phi$-uncolored edge and let $i$ be the least index such that $e = e_i$. If $i > 0$, then $e$ is still uncolored by $\Shift_{i-1}(\phi, C)$, which means that the pair $(e_{i-1}, e_i)$ is not $\Shift_{i-1}(\phi, C)$-shiftable; a contradiction. Part \ref{item:end} follows by applying \ref{item:start} with $\Shift(\phi, C)$ and $C^\ast$ in place of $\phi$ and $C$.
	\end{scproof}
	
	We say that a chain $C$ is \emphd{$\phi$-happy} if $C$ is $\phi$-shiftable and the edge $\eend(C)$ is $\Shift(\phi, C)$-happy.
	
	\begin{fact}\label{fact:happy}
		Let $\phi$ be a proper partial coloring and let $C$ be a chain. If $C$ is $\phi$-happy, then the graph $H(C)$ is augmenting for $\phi$.
	\end{fact}
	\begin{scproof}
		The coloring $\phi$ can be modified to include all the edges in $E(C)$ as follows: First, replace $\phi$ by $\psi \defeq \Shift(\phi, C)$. Now the only uncolored edge in $E(C)$ is $\eend(C)$, which, by assumption, is $\psi$-happy, so $\psi$ can be extended to include it.
	\end{scproof}
	
	Thus, our goal becomes to develop a technique for constructing ``short'' happy chains.
	
	\subsection{Initial segments and combinations}
	
	For a chain $C = (e_0, \ldots, e_{\ell-1})$ and $1 \leq k \leq \ell$, the \emphd{initial segment} of $C$ of length $k$ is the chain
	\[
	\rest{C}{k} \,\defeq\, (e_0, \ldots, e_{k-1}).
	\]
	Note that if $C$ is $\phi$-shiftable for some proper partial coloring $\phi$, then every initial segment of $C$ is also $\phi$-shiftable. Observe also that $\Shift(\phi, \rest{C}{k}) = \Shift_{k-1}(\phi, C)$ for all $1 \leq k \leq \ell$. We additionally consider an operation that is in some sense the opposite of taking initial segments. Namely, let $C_0$ and $C_1$ be chains such that $\eend(C_0) = \start(C_1)$. Then we can \emphd{combine} $C_0$ and $C_1$ into a single chain $C_0 + C_1$ by identifying the last edge of $C_0$ with the first edge of $C_1$. More precisely, if $C_0 = (e_0, \ldots, e_{\ell-1})$, $C_1 = (h_0, \ldots, h_{k-1})$, and $e_{\ell-1} = h_0$, then we set
	\[
	C_0 + C_1 \,\defeq\, (e_0, \ldots, e_{\ell-1}, h_1, \ldots, h_{k-1}).
	\]
	Notice that $C_0$ is an initial segment of $C_0 + C_1$. It is clear that for any proper partial coloring $\phi$, we have
	\[
	\Shift(\phi, C_0 + C_1) \,=\, \Shift(\Shift(\phi, C_0), C_1),
	\]
	and if $C_0$ is $\phi$-shiftable and $C_1$ is $\Shift(\phi, C_0)$-shiftable, then $C_0 + C_1$ is also $\phi$-shiftable. 

	\subsection{Path-chains}\label{subsec:path}
	
	A particular type of a chain is a \emph{path-chain}, i.e., a nonempty sequence of edges that forms a path in $G$ in the usual graph-theoretic sense. Explicitly, a chain $P = (e_0, \ldots, e_{\ell - 1})$ is called a \emphd{path-chain}, or simply a \emphd{path}, if there is a sequence $(x_0, x_1, \ldots, x_{\ell})$ of pairwise distinct vertices such that $e_i = x_i x_{i+1}$ for all $1 \leq i \leq \ell - 1$. Such a sequence $(x_0, x_1, \ldots, x_{\ell})$ is necessarily unique unless $\ell = 1$, while if $P = (xy)$ is a path-chain of length $1$, then there are two such sequences, namely $(x,y)$ and $(y,x)$. When dealing with such length-$1$ path-chains, we shall still tacitly assume that one of these two sequences is chosen, so that every path-chain $P$ has a well-defined first and a well-defined last vertex, denoted by $\vstart(P)$ and $\vend(P)$, respectively. 
	
	The aim of this subsection is to describe a particular way of constructing shiftable path-chains by considering edges of two specific colors. For the remainder of \S\ref{subsec:path}, we fix the following data:
	\begin{itemize}
		\item a proper partial coloring $\phi$;
		
		\item a pair of distinct colors $\alpha$, $\beta \in [\Delta + 1]$. 
	\end{itemize}
	Let $G(\phi, \alpha \beta)$ denote the spanning subgraph of $G$ with
	\[
	E(G(\phi, \alpha \beta)) \,\defeq\, \set{e \in E \,:\, \phi(e) \in \set{\alpha, \beta}}.
	\]
	Since $\phi$ is proper, the maximum degree of $G(\phi, \alpha \beta)$ is at most $2$, which means that every connected component of $G(\phi, \alpha \beta)$ is either a cycle or a path \ep{where we view isolated vertices as paths of length $0$}. For $x \in V$, let $G(x, \phi, \alpha \beta)$ denote the connected component of $x$ in $G(\phi, \alpha \beta)$, and let $\deg(x, \phi, \alpha \beta)$ be the degree of $x$ in $G(\phi,\alpha \beta)$. Note that $\deg(x, \phi,\alpha \beta) < 2$ if and only if at least one of $\alpha$, $\beta$ is missing at $x$. We say that two vertices $x$, $y \in V$ are \emphd{$(\phi, \alpha \beta)$-related} if $G(x, \phi,\alpha \beta) = G(y, \phi,\alpha \beta)$.
	
\begin{figure}[h]
		\centering
		\begin{subfigure}[t]{.4\textwidth}
			\centering
			\begin{tikzpicture}
			\node[draw=none,minimum size=2.5cm,regular polygon,regular polygon sides=7] (P) {};

			\node[circle,fill=black,draw,inner sep=0pt,minimum size=4pt] (x) at (P.corner 4) {};
			\node[circle,fill=black,draw,inner sep=0pt,minimum size=4pt] (y) at (P.corner 5) {};
			\node[circle,fill=black,draw,inner sep=0pt,minimum size=4pt] (a) at (P.corner 6) {};
			\node[circle,fill=black,draw,inner sep=0pt,minimum size=4pt] (b) at (P.corner 7) {};
			\node[circle,fill=black,draw,inner sep=0pt,minimum size=4pt] (c) at (P.corner 1) {};
			\node[circle,fill=black,draw,inner sep=0pt,minimum size=4pt] (d) at (P.corner 2) {};
			\node[circle,fill=black,draw,inner sep=0pt,minimum size=4pt] (e) at (P.corner 3) {};
			
			\node[anchor=north] at (x) {$x$};
			\node[anchor=north] at (y) {$y$};
			
			\draw[thick,dotted] (x) to (y);
			\draw[thick] (y) to node[midway,inner sep=1pt,outer sep=1pt,minimum size=4pt,fill=white] {$\alpha$} (a);
			\draw[thick] (a) to node[midway,inner sep=1pt,outer sep=1pt,minimum size=4pt,fill=white] {$\beta$} (b);
			\draw[thick] (b) to node[midway,inner sep=1pt,outer sep=1pt,minimum size=4pt,fill=white] {$\alpha$} (c);
			\draw[thick] (c) to node[midway,inner sep=1pt,outer sep=1pt,minimum size=4pt,fill=white] {$\beta$} (d);
			\draw[thick] (d) to node[midway,inner sep=1pt,outer sep=1pt,minimum size=4pt,fill=white] {$\alpha$} (e);
			\draw[thick] (e) to node[midway,inner sep=1pt,outer sep=1pt,minimum size=4pt,fill=white] {$\beta$} (x);
			\end{tikzpicture}
			\caption{The edge $xy$ is not $(\phi, \alpha \beta)$-successful.}
		\end{subfigure}%
		\qquad%
		\begin{subfigure}[t]{.4\textwidth}
			\centering
			\begin{tikzpicture}
			\node[draw=none,minimum size=2.5cm,regular polygon,regular polygon sides=7] (P) {};

			\node[circle,fill=black,draw,inner sep=0pt,minimum size=4pt] (x) at (P.corner 4) {};
			\node[circle,fill=black,draw,inner sep=0pt,minimum size=4pt] (y) at (P.corner 5) {};
			\node[circle,fill=black,draw,inner sep=0pt,minimum size=4pt] (a) at (P.corner 6) {};
			\node[circle,fill=black,draw,inner sep=0pt,minimum size=4pt] (b) at (P.corner 7) {};
			\node[circle,fill=black,draw,inner sep=0pt,minimum size=4pt] (c) at (P.corner 1) {};
			\node[circle,fill=black,draw,inner sep=0pt,minimum size=4pt] (d) at (P.corner 2) {};
			\node[circle,fill=black,draw,inner sep=0pt,minimum size=4pt] (e) at (P.corner 3) {};
			\node[circle,fill=black,draw,inner sep=0pt,minimum size=4pt] (f) at (-2.2,0) {}; 
			
			\node[anchor=north] at (x) {$x$};
			\node[anchor=north] at (y) {$y$};
			
			\draw[thick,dotted] (x) to (y);
			\draw[thick] (y) to node[midway,inner sep=1pt,outer sep=1pt,minimum size=4pt,fill=white] {$\alpha$} (a);
			\draw[thick] (a) to node[midway,inner sep=1pt,outer sep=1pt,minimum size=4pt,fill=white] {$\beta$} (b);
			\draw[thick] (b) to node[midway,inner sep=1pt,outer sep=1pt,minimum size=4pt,fill=white] {$\alpha$} (c);
			\draw[thick] (c) to node[midway,inner sep=1pt,outer sep=1pt,minimum size=4pt,fill=white] {$\beta$} (d);
			\draw[thick] (d) to node[midway,inner sep=1pt,outer sep=1pt,minimum size=4pt,fill=white] {$\alpha$} (e);
			\draw[thick] (e) to node[midway,inner sep=1pt,outer sep=1pt,minimum size=4pt,fill=white] {$\beta$} (f);
			\end{tikzpicture}
			\caption{The edge $xy$ is $(\phi, \alpha \beta)$-successful.}
		\end{subfigure}
		\caption{The chain $P(x, y, \phi, \alpha\beta)$.}\label{fig:path}
	\end{figure}
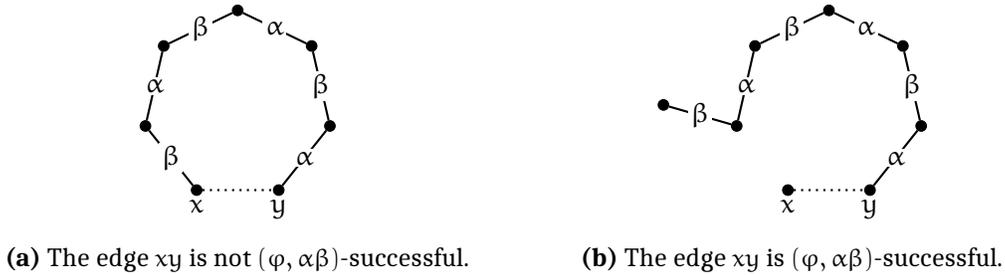

	Suppose that $xy \in E$ is an uncolored edge that is not $\phi$-happy. We say that $xy$ is:
	\begin{itemize}
		\item \emphd{$(\phi, \alpha \beta)$-hopeful} if $\deg(x, \phi, \alpha \beta) < 2$ and $\deg(y, \phi, \alpha \beta) < 2$; 
		
		\item \emphd{$(\phi, \alpha \beta)$-successful} if it is $(\phi, \alpha \beta)$-hopeful and $x$ and $y$ are not $(\phi, \alpha \beta)$-related.
	\end{itemize}
	For a $(\phi, \alpha \beta)$-hopeful edge $xy$, we define a chain $P(x, y, \phi, \alpha \beta)$ as follows. Since 
	$\deg(y, \phi, \alpha \beta) < 2$, 
	the graph $G(y, \phi, \alpha\beta)$ is a path, one of whose endpoints is $y$. Let $(e_1, \ldots, e_k)$ be the sequence in which the edges of this path appear as it is traversed starting from $y$. Note that $k \geq 1$, since otherwise $\set{\alpha, \beta} \subseteq M(\phi, y)$ and $\set{\alpha, \beta} \cap M(\phi, x) \neq \0$, which would mean that the edge $xy$ is $\phi$-happy. The sequence $(e_1, \ldots, e_k)$ is \emph{$\alpha\beta$-alternating}, in the sense that the sequence of colors $\phi(e_1)$, $\phi(e_2)$, $\phi(e_3)$, $\phi(e_4)$, \ldots{} has the form $\alpha$, $\beta$, $\alpha$, $\beta$, \ldots{} or $\beta$, $\alpha$, $\beta$, $\alpha$, \ldots. 
	Then we define
	\[
	P(x, y, \phi, \alpha\beta) \,\defeq\, (xy, e_1, \ldots, e_k).
	\]
	This construction is illustrated in Fig.~\ref{fig:path}.
	
	\begin{fact}\label{fact:hope}
		If $xy$ is a $(\phi, \alpha \beta)$-hopeful edge, then the chain $P(x, y, \phi, \alpha\beta)$ is $\phi$-shiftable.
	\end{fact}
	\begin{scproof}
		Since the edge $xy$ is not $\phi$-happy, $M(\phi, x) \cap M(\phi, y) = \0$. As $xy$ is $(\phi, \alpha \beta)$-hopeful, this implies that $x$ misses one of the colors $\alpha$, $\beta$ and $y$ misses the other one. Without loss of generality, assume that $\alpha \in M(\phi, x)$ and $\beta \in M(\phi, y)$. Set $e_0 \defeq xy$ and $P \defeq P(x, y, \phi, \alpha\beta) = (e_0, e_1, \ldots, e_k)$. Then $e_0$ is uncolored and $\phi(e_1) = \alpha$, $\phi(e_2) = \beta$, $\phi(e_3) = \alpha$, $\phi(e_4) = \beta$, etc. Since $\alpha \in M(\phi, x)$, the pair of edges $(e_0, e_1)$ is $\phi$-shiftable. Straightforward induction then shows that the whole chain $P$ is $\phi$-shiftable.
	\end{scproof}
	
	If an edge $xy$ is not just $(\phi, \alpha \beta)$-hopeful, but in fact $(\phi, \alpha \beta)$-successful, then $x$ cannot appear on the path $G(y, \phi, \alpha\beta)$, which implies that $P(x, y, \phi, \alpha\beta)$ is a path-chain. Note that the first vertex of $P(x, y, \phi, \alpha\beta)$ is $x$, while the last one is the endpoint of the path $G(y, \phi, \alpha\beta)$ distinct from $y$.
	
	\begin{fact}\label{fact:happy_path}
		If $xy$ is a $(\phi, \alpha \beta)$-successful edge, then the path-chain $P(x, y, \phi, \alpha\beta)$ is $\phi$-happy.
	\end{fact}
	\begin{scproof}
		Set $e_0 \defeq xy$ and $P \defeq P(x, y, \phi, \alpha\beta) = (e_0, e_1, \ldots, e_k)$. We already know, from Fact~\ref{fact:hope}, that $P$ is $\phi$-shiftable. In the coloring $\Shift(\phi, P)$, both endpoints of the edge $e_k$ miss the same one of the colors $\alpha$, $\beta$, so $e_k$ is $\Shift(\phi, P)$-happy, as desired.
	\end{scproof}
	
	\subsection{Fans}
	
	Another useful type of a chain is a \emph{fan}. A \emphd{fan} is a chain of the form $F = (xy_0, xy_1, \ldots, x y_{\ell -1})$, where $x \in V$ and $y_0$, \ldots, $y_{\ell - 1}$ are pairwise distinct neighbors of $x$. Given a fan $F = (xy_0, xy_1, \ldots, x y_{\ell -1})$, we call the vertex $x$ the \emphd{pivot} of $F$ and write $\pivot(F) \defeq x$. We also refer to $y_{\ell - 1}$ as the last vertex of $F$ and write $\vend(F) \defeq y_{\ell - 1}$. As in the case of path-chains, this creates some ambiguity when $F = (xy)$ is a fan of length $1$. When dealing with such a length-$1$ fan, we shall still tacitly assume that it has been decided which one of $x$, $y$ is the pivot. Note that for every fan $F$, its last edge $\eend(F)$ is $\set{\pivot(F), \vend(F)}$.
	
	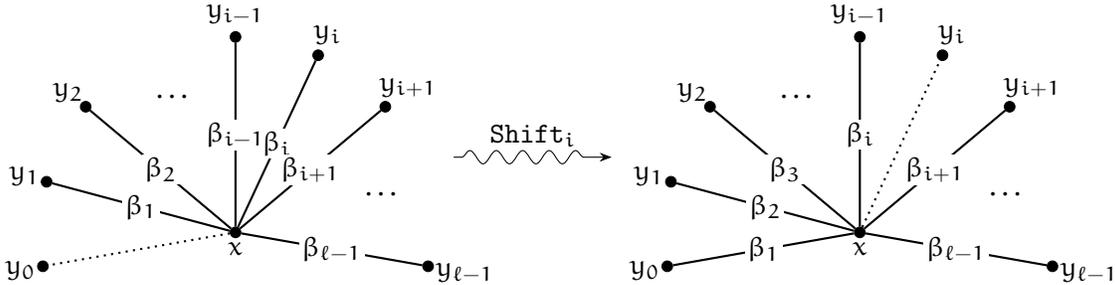
\begin{figure}[h]
		\centering
		\begin{tikzpicture}
		\begin{scope}
		\node[circle,fill=black,draw,inner sep=0pt,minimum size=4pt] (x) at (0,0) {};
		\node[anchor=north] at (x) {$x$};
		
		\coordinate (O) at (0,0);
		\def\radius{2.6cm}
		
		\node[circle,fill=black,draw,inner sep=0pt,minimum size=4pt] (y0) at (190:\radius) {};
		\node at (190:2.9) {$y_0$};
		
		\node[circle,fill=black,draw,inner sep=0pt,minimum size=4pt] (y1) at (165:\radius) {};
		\node at (165:2.9) {$y_1$};
		
		\node[circle,fill=black,draw,inner sep=0pt,minimum size=4pt] (y2) at (140:\radius) {};
		\node at (140:2.9) {$y_2$};
		
		\node[circle,fill=black,draw,inner sep=0pt,minimum size=4pt] (y4) at (90:\radius) {};
		\node at (90:2.9) {$y_{i-1}$};
		
		\node[circle,fill=black,draw,inner sep=0pt,minimum size=4pt] (y5) at (65:\radius) {};
		\node at (65:2.9) {$y_i$};
		
		\node[circle,fill=black,draw,inner sep=0pt,minimum size=4pt] (y6) at (40:\radius) {};
		\node at (40:3) {$y_{i+1}$};
		
		\node[circle,fill=black,draw,inner sep=0pt,minimum size=4pt] (y8) at (-10:\radius) {};
		\node at (-10:3.1) {$y_{\ell-1}$};
		
		\node[circle,inner sep=0pt,minimum size=4pt] at (115:2) {$\ldots$}; 
		\node[circle,inner sep=0pt,minimum size=4pt] at (15:2) {$\ldots$}; 
		
		\draw[thick,dotted] (x) to (y0);
		\draw[thick] (x) to node[midway,inner sep=1pt,outer sep=1pt,minimum size=4pt,fill=white] {$\beta_1$} (y1);
		\draw[thick] (x) to node[midway,inner sep=1pt,outer sep=1pt,minimum size=4pt,fill=white] {$\beta_2$} (y2);
		
		\draw[thick] (x) to node[midway,inner sep=1pt,outer sep=1pt,minimum size=4pt,fill=white] {$\beta_{i-1}$} (y4);
		\draw[thick] (x) to node[midway,inner sep=1pt,outer sep=1pt,minimum size=4pt,fill=white] {$\beta_{i}$} (y5);
		\draw[thick] (x) to node[midway,inner sep=1pt,outer sep=1pt,minimum size=4pt,fill=white] {$\beta_{i+1}$} (y6);
		
		\draw[thick] (x) to node[midway,inner sep=1pt,outer sep=1pt,minimum size=4pt,fill=white] {$\beta_{\ell - 1}$} (y8);
		\end{scope}
		
		\draw[-{Stealth[length=1.6mm]},decoration = {snake,pre length=3pt,post length=7pt,},decorate] (2.9,1) to node[midway,anchor=south]{$\Shift_i$} (5,1);
		
		\begin{scope}[xshift=8.3cm]
		\node[circle,fill=black,draw,inner sep=0pt,minimum size=4pt] (x) at (0,0) {};
		\node[anchor=north] at (x) {$x$};
		
		\coordinate (O) at (0,0);
		\def\radius{2.6cm}
		
		\node[circle,fill=black,draw,inner sep=0pt,minimum size=4pt] (y0) at (190:\radius) {};
		\node at (190:2.9) {$y_0$};
		
		\node[circle,fill=black,draw,inner sep=0pt,minimum size=4pt] (y1) at (165:\radius) {};
		\node at (165:2.9) {$y_1$};
		
		\node[circle,fill=black,draw,inner sep=0pt,minimum size=4pt] (y2) at (140:\radius) {};
		\node at (140:2.9) {$y_2$};
		
		\node[circle,fill=black,draw,inner sep=0pt,minimum size=4pt] (y4) at (90:\radius) {};
		\node at (90:2.9) {$y_{i-1}$};
		
		\node[circle,fill=black,draw,inner sep=0pt,minimum size=4pt] (y5) at (65:\radius) {};
		\node at (65:2.9) {$y_i$};
		
		\node[circle,fill=black,draw,inner sep=0pt,minimum size=4pt] (y6) at (40:\radius) {};
		\node at (40:3) {$y_{i+1}$};
		
		\node[circle,fill=black,draw,inner sep=0pt,minimum size=4pt] (y8) at (-10:\radius) {};
		\node at (-10:3.1) {$y_{\ell-1}$};
		
		\node[circle,inner sep=0pt,minimum size=4pt] at (115:2) {$\ldots$}; 
		\node[circle,inner sep=0pt,minimum size=4pt] at (15:2) {$\ldots$}; 
		
		\draw[thick] (x) to node[midway,inner sep=1pt,outer sep=1pt,minimum size=4pt,fill=white] {$\beta_1$} (y0);
		\draw[thick] (x) to node[midway,inner sep=1pt,outer sep=1pt,minimum size=4pt,fill=white] {$\beta_2$} (y1);
		\draw[thick] (x) to node[midway,inner sep=1pt,outer sep=1pt,minimum size=4pt,fill=white] {$\beta_3$} (y2);
		
		\draw[thick] (x) to node[midway,inner sep=1pt,outer sep=1pt,minimum size=4pt,fill=white] {$\beta_{i}$} (y4);
		\draw[thick,dotted] (x) to (y5);
		\draw[thick] (x) to node[midway,inner sep=1pt,outer sep=1pt,minimum size=4pt,fill=white] {$\beta_{i+1}$} (y6);
		
		\draw[thick] (x) to node[midway,inner sep=1pt,outer sep=1pt,minimum size=4pt,fill=white] {$\beta_{\ell - 1}$} (y8);
		\end{scope}
		\end{tikzpicture}
		\caption{The process of shifting a fan for $i$ steps.}\label{fig:fan}
	\end{figure}
	
	\vspace*{-5pt}
	
	\begin{fact}
		Let $\phi$ be a proper partial coloring and consider a fan $F = (xy_0, xy_1, \ldots, x y_{\ell -1})$. Then $F$ is $\phi$-shiftable if and only if the following two conditions are satisfied:
		\begin{enumerate}[label=\ep{\normalfont\roman*}]
			\item\label{item:fan_i} $xy_0 \not \in \dom(\phi)$, while $xy_i \in \dom(\phi)$ for all $1 \leq i \leq \ell - 1$;
			
			\item\label{item:fan_ii} for each $0 \leq i \leq \ell - 2$, $\phi(xy_{i+1}) \in M(\phi, y_i)$.
		\end{enumerate}
	\end{fact}
	\begin{scproof}
		Condition \ref{item:fan_i} is necessary for $F$ to be $\phi$-shiftable due to Fact~\ref{fact:shift}, so we may assume that \ref{item:fan_i} holds. For brevity, set $\psi_i \defeq \Shift_i(\phi, F)$ for each $0 \leq i \leq \ell -1$. By definition, $F$ is $\phi$-shiftable if for all $0 \leq i \leq \ell - 2$, the pair $(xy_i, xy_{i+1})$ is $\psi_i$-shiftable. It is straightforward to check that $\psi_i(xy_i) = \blank$, while $\psi_i(xy_{i+1}) = \phi(xy_{i+1})$ \ep{see Fig.~\ref{fig:fan}}. Thus, $(xy_i, xy_{i+1})$ is $\psi_i$-shiftable if and only if $\phi(xy_{i+1}) \in M(\psi_i, y_i)$. For $i = 0$, this is just saying that $\phi(xy_1) \in M(\phi, y_0)$, as desired. For $i > 0$, observe that $M(\psi_i, y_i) = M(\phi, y_i) \cup \set{\phi(xy_i)}$. Since the coloring $\phi$ is proper, $\phi(xy_{i+1}) \neq \phi(xy_{i})$, so $\phi(xy_{i+1}) \in M(\psi_i, y_i)$ if and only if $\phi(xy_{i+1}) \in M(\phi, y_i)$, and we are done.
	\end{scproof}

	Let $\phi$ be a partial proper coloring and let $\alpha$, $\beta \in [\Delta + 1]$ be two distinct colors.
	Let $F$ be a $\phi$-shiftable fan that is not $\phi$-happy. Set $x \defeq \pivot(F)$ and $y \defeq \vend(F)$ \ep{so $\eend(F) = xy$}. By analogy with the definitions in \S\ref{subsec:path}, we say that $F$ is:
	\begin{itemize}
		\item \emphd{$(\phi, \alpha \beta)$-hopeful} if $\deg(x, \phi, \alpha \beta) < 2$ and $\deg(y, \phi, \alpha \beta) < 2$;
		
		\item \emphd{$(\phi, \alpha \beta)$-successful} if it is $(\phi, \alpha \beta)$-hopeful and $x$ and $y$ are not $(\Shift(\phi, F), \alpha \beta)$-related.
	\end{itemize}
	
	\begin{fact}\label{fact:prom}
		Let $F$ be a fan and set $x \defeq \pivot(F)$ and $y \defeq \vend(F)$. Define $\psi \defeq \Shift(\phi, F)$. If $F$ is $(\phi, \alpha \beta)$-hopeful \ep{resp.~$(\phi, \alpha \beta)$-successful}, then the edge $xy$ is $(\psi, \alpha \beta)$-hopeful \ep{resp.~$(\psi, \alpha\beta)$-successful}.
	\end{fact}
	\begin{scproof}
		Note that $M(\psi, x) = M(\phi, x)$ and $M(\psi, y) \supseteq M(\phi, y)$ \ep{the inclusion is, in fact, strict unless $F$ is of length $1$}. This implies that if $F$ is $(\phi, \alpha \beta)$-hopeful, then $\deg(x, \psi, \alpha \beta) = \deg(x, \phi, \alpha \beta) < 2$ and $\deg(y, \psi, \alpha \beta) \leq \deg(y, \phi, \alpha \beta) < 2$, and hence the edge $xy$ is $(\psi, \alpha \beta)$-hopeful, as desired. The statement for $(\phi, \alpha \beta)$-successful fans  follows.
	\end{scproof}
	
	Next we establish two lemmas that allow us to find successful or, at least, hopeful fans. Modulo our choice of terminology, the first of these lemmas was already known to Vizing \cite{Viz64} \ep{see \cite{Die} for a textbook presentation}, while the second one comes from the work of Greb\'\i{}k and Pikhurko \cite{GP}.
	
	\begin{lemma}[\textls{First fan lemma}]\label{lemma:simple_fan}
		Let $\phi$ be a partial proper coloring and let $xy \in E \setminus \dom(\phi)$ be an uncolored edge. Then there exists a $\phi$-shiftable fan $F$ with pivot $x$ and $\start(F) = xy$ such that:
		\begin{itemize}
			\item either $F$ is $\phi$-happy,
			
			\item or, for some distinct $\alpha$, $\beta \in [\Delta + 1]$, $F$ is $(\phi, \alpha \beta)$-successful.
		\end{itemize}
	\end{lemma}
	\begin{scproof}
		Assume that there is no $\phi$-happy fan with pivot $x$ and starting edge $xy$. For each neighbor $z$ of $x$, fix an arbitrary color $\beta(z) \in M(\phi, z)$. We shall recursively build a sequence $y_0$, $y_1$, \ldots{} of pairwise distinct neighbors of $x$ such that the fan $(xy_0, xy_1, \ldots)$ is $\phi$-shiftable. Start by setting $y_0 \defeq y$. Once the vertices $y_0$, \ldots, $y_i$ have been determined, consider the color $\beta(y_i)$. If $\beta(y_i) \in M(\phi, x)$, then the fan $(xy_0, \ldots, xy_i)$ is $\phi$-happy, which contradicts our assumption. Thus, $\beta(y_i) \not \in M(\phi, x)$, which means that $x$ has a \ep{unique} neighbor $z$ with $\phi(xz) = \beta(y_i)$. Note that $z \neq y_i$, because the color $\beta(y_i)$ is missing at $y_i$ but not at $z$. Similarly, $z \neq y_0$ because the edge $xy_0$ is uncolored, while $xz$ is not. Now, if $z \in \set{y_1, \ldots, y_{i-1}}$, then we stop the construction, and otherwise we set $y_{i+1} \defeq z$ and proceed to the next step. Since $x$ has only finitely many neighbors, this construction must terminate. That is, we shall eventually obtain a sequence $y_0$, \ldots, $y_k$ of distinct neighbors of $x$ such that:
		\begin{enumerate}[label=\ep{\normalfont\roman*}]
			\item $y_0 = y$;
			\item\label{item:not_missing_at_x} for all $0 \leq i \leq k$, $\beta(y_i) \not \in M(\phi, x)$;
			\item\label{item:step} for all $0 \leq i \leq k - 1$, $\phi(xy_{i+1}) = \beta(y_i)$;
			\item\label{item:back} $\beta(y_k) = \phi(xy_j)$ for some $1 \leq j \leq k-1$.
		\end{enumerate}
		Set $\beta \defeq \beta(y_k)$. If $j$ is the index from~\ref{item:back}, then \ref{item:step} yields $\beta = \phi(xy_j) = \beta(y_{j-1})$; that is, $\beta$ is missing both at $y_k$ and at $y_{j-1}$. Now let $\alpha$ be any color in $M(\phi, x)$. By \ref{item:not_missing_at_x}, $\alpha \neq \beta$. Define
		\[F \,\defeq\, (xy_0, \ldots, xy_k) \qquad \text{and} \qquad F' \,\defeq\, (xy_0, \ldots, xy_{j-1}) \,=\, \rest{F}{j}.\]
		By construction, both $F$ and $F'$ are $(\phi, \alpha \beta)$-hopeful fans. We claim that at lest one of $F$, $F'$ is $(\phi, \alpha \beta)$-successful. For brevity, let
		\[\psi \,\defeq\, \Shift(\phi, F) \qquad \text{and} \qquad \psi' \,\defeq\, \Shift(\phi, F') \,=\, \Shift_{j-1}(\phi, F).\]
		Suppose that $F$ is not $(\phi, \alpha \beta)$-successful. This means that the vertices $x$ and $y_k$ are $(\psi, \alpha \beta)$-related. By Fact~\ref{fact:prom}, the edge $xy_k$ is $(\psi, \alpha \beta)$-hopeful, 
		so $G_{\alpha \beta}(\psi, x)$ must be a path with endpoints $x$ and $y_k$. Furthermore, $\psi(xy_{j-1}) = \beta$, so the \ep{unique} edge of this path incident to $x$ is $xy_{j-1}$. All the other edges of this path avoid the vertex $x$, which means that their colors cannot be changed by shifting a fan with pivot $x$. In particular, they have the same colors in $\phi$, $\psi$, and $\psi'$, which implies that the vertices $y_{j-1}$ and $y_k$ are $(\psi', \alpha \beta)$-related. Since both these vertices miss the color $\beta$ in $\psi'$, the graph $G_{\alpha \beta}(\psi', y_{j-1})$ is a path with endpoints $y_{j-1}$ and $y_k$. This path cannot include $x$ \ep{because $x$ is missing $\alpha$}, and thus $y_{j-1}$ and $x$ are not $(\psi', \alpha \beta)$-related, as desired.
	\end{scproof}
	
	\begin{lemma}[\textls{Second fan lemma}]\label{lemma:fancy_fan}
		Let $\phi$ be a partial proper coloring and let $xy \in E \setminus \dom(\phi)$ be an uncolored edge. Fix a pair of distinct colors $\alpha$, $\beta \in [\Delta + 1]$ such that $\deg(x, \phi, \alpha\beta) = 1$. 
		Then there exists a $\phi$-shiftable fan $F$ with pivot $x$ and $\start(F) = xy$ such that:
		\begin{itemize}
			\item either $F$ is $\phi$-happy,
			
			\item or, for some distinct $\gamma$, $\delta \in [\Delta + 1]$ with $\set{\gamma, \delta} \cap \set{\alpha, \beta} = \0$, $F$ is $(\phi, \gamma \delta)$-successful,
			
			\item or $F$ is $(\phi, \alpha \beta)$-hopeful, $\vend(F) \neq y$, and no edge in $F$ is colored $\alpha$ or $\beta$ by $\phi$.
		\end{itemize}
	\end{lemma}
	\begin{scproof}
		Again, we may assume that there is no $\phi$-happy fan with pivot $x$ and starting edge $xy$. Since $\deg(x, \phi, \alpha\beta) = 1$, precisely one of the colors $\alpha$, $\beta$ is missing at $x$. For concreteness, suppose that $\alpha \in M(\phi, x)$ and $\beta \not \in M(\phi, x)$. For each neighbor $z$ of $x$, fix a color $\delta(z) \in M(\phi, z)$; it does not matter which color from $M(\phi, z)$ we pick, except that we make sure that $\delta(y) \neq \beta$ \ep{this is possible since $y$ is incident to an uncolored edge and thus $|M(\phi, y)| \geq 2$}.  As in the proof of Lemma~\ref{lemma:simple_fan}, we build a sequence $y_0$, $y_1$, \ldots{} of pairwise distinct neighbors of $x$ such that the fan $(xy_0, xy_1, \ldots)$ is $\phi$-shiftable. Start by setting $y_0 \defeq y$. Once the vertices $y_0$, \ldots, $y_i$ have been determined, consider the color $\delta(y_i)$. If $\delta(y_i) \in M(\phi, x)$, then the fan $(xy_0, \ldots, xy_i)$ is $\phi$-happy, which contradicts our assumption. Thus, $\delta(y_i) \not \in M(\phi, x)$. At this point, if $\delta(y_i) = \beta$, then we stop the construction. Otherwise we proceed as in the proof of Lemma~\ref{lemma:simple_fan}. In other words, if $\delta(y_i) \neq \beta$, then we let $z$ be the \ep{unique} neighbor of $x$ such that $\phi(xz) = \delta(y_i)$. If $z \in \set{y_1, \ldots, y_{i-1}}$, then we stop, and otherwise, we let $y_{i+1} \defeq z$ and move on to the next step. Since $x$ has only finitely many neighbors, this construction must terminate. Hence, we eventually obtain a sequence $y_0$, \ldots, $y_k$ of distinct neighbors of $x$ such that:
		\begin{enumerate}[label=\ep{\normalfont\roman*}]
			\item $y_0 = y$;
			\item\label{item:not_missing_at_x_1} for all $0 \leq i \leq k$, $\delta(y_i) \not \in M(\phi, x)$;
			\item for all $0 \leq i \leq k - 1$, $\phi(xy_{i+1}) = \delta(y_i)$;
			\item\label{item:not_beta} for all $0 \leq i \leq k - 1$, $\delta(y_i) \neq \beta$;
			\item\label{item:back_1} either $\delta(y_k) = \beta$, or else, $\delta(y_k) = \phi(xy_j)$ for some $1 \leq j \leq k-1$.
		\end{enumerate}
		Now we consider two cases.
		
		\setcounter{case}{0}
		\begin{case}
			$\delta(y_k) = \beta$.
		\end{case}
		
		Then the fan $F \defeq (xy_0, \ldots, xy_k)$ is, by construction, $(\phi, \alpha \beta)$-hopeful. Furthermore, since $\delta(y) \neq \beta$, we have $\vend(F) = y_k \neq y$. Finally, for each $1 \leq i \leq k$, we have $\phi(xy_i) = \delta(y_{i-1})$, and since $\delta(y_{i-1}) \neq \beta$ by \ref{item:not_beta} and $\delta(y_{i-1}) \not \in M(\phi, x)$  by \ref{item:not_missing_at_x_1}, we conclude that no edge in $F$ is colored $\alpha$ or $\beta$.
		
		\begin{case}
			$\delta(y_k) \neq \beta$ and $\delta(y_k) = \phi(xy_j)$ for some $1 \leq j \leq k-1$.
		\end{case}
		
		Set $\delta \defeq \delta(y_k)$ and let $j$ be the index from \ref{item:back_1}, so that the color $\delta$ is missing both at $y_k$ and at $y_{j-1}$. Since $x$ is incident to an uncolored edge, we can pick a color $\gamma \in M(x, \phi)$ distinct from $\alpha$. Now $\set{\gamma, \delta} \cap \set{\alpha, \beta} = \0$ and the argument from the proof of Lemma~\ref{lemma:simple_fan} shows that at least one of the fans
		\[F \,\defeq\, (xy_0, \ldots, xy_k) \qquad \text{and} \qquad F' \,\defeq\, (xy_0, \ldots, xy_{j-1})\]
		is $(\phi, \gamma\delta)$-successful, as desired.
	\end{scproof}
	
	As an illustration, we can now use Lemma~\ref{lemma:simple_fan} to complete the proof of Vizing's theorem:
	
	\begin{corl}[{to Lemma~\ref{lemma:simple_fan}}]\label{corl:vizing}
		For every proper partial coloring $\phi$ and every uncolored edge $e = xy$, there is a connected augmenting subgraph $H \subseteq G$ whose only uncolored edge is $e$. 
	\end{corl}
	\begin{scproof}
		By Fact~\ref{fact:happy}, it suffices to find a $\phi$-happy chain $C$ with $\start(C) = e$. To this end, let $F$ be a fan satisfying the conclusion of Lemma~\ref{lemma:simple_fan}. If $F$ is $\phi$-happy, then we are done, so assume that instead there are distinct colors $\alpha$, $\beta \in [\Delta + 1]$ such that $F$ is $(\phi, \alpha \beta)$-successful. Let $z \defeq \vend(F)$ and set $\psi \defeq \Shift(\phi, F)$. By Fact~\ref{fact:prom}, the edge $xz$ is $(\psi, \alpha \beta)$-successful, so, by Fact~\ref{fact:happy_path}, the path-chain $P \defeq P(\psi, x, z, \alpha \beta)$ is $\psi$-happy. Therefore, the chain $F + P$ is as desired.
	\end{scproof}
	

	\section{The Multi-Step Vizing Algorithm}\label{sec:const}
	
	\subsection{The sets $R^\shortrightarrow(x, \phi)$ and $R^\shortleftarrow(y, \phi)$}
	
	Before we begin describing our construction of multi-step Vizing chains, we need to introduce some notation that will be useful for keeping track of intersections between different alternating paths. Let $\phi$ be a proper partial coloring and let $\alpha$, $\beta \in [\Delta + 1]$ be distinct colors. Given $x \in V$, we let $R^\shortrightarrow(x, \phi, \alpha\beta) \subseteq V$ denote the set of all vertices $y$ that satisfy at least one of the following two conditions:
	\begin{itemize}
		\item $y = x$ or $y$ is adjacent to $x$,
		
		\item or there is a neighbor $z$ of $x$ such that $\deg(z, \phi, \alpha\beta) < 2$ and $z$ and $y$ are $(\phi, \alpha\beta)$-related.
	\end{itemize}
	We also define $R^\shortleftarrow(y, \phi, \alpha\beta) \defeq \set{x \in V \,:\, y \in R^\shortrightarrow(x, \phi, \alpha\beta)}$ and set
	\[
	R^\shortrightarrow(x, \phi) \,\defeq\, \bigcup_{\alpha \neq \beta} R^\shortrightarrow(x, \phi, \alpha\beta) \quad \text{and} \quad R^\shortleftarrow(y, \phi) \,\defeq\, \bigcup_{\alpha \neq \beta} R^\shortleftarrow(y, \phi, \alpha\beta).
	\]
	
	\begin{lemma}\label{lemma:indegree}
		For each proper partial coloring $\phi$ and every vertex $y \in V$, we have $|R^\shortleftarrow(y, \phi)| \leq (\Delta + 1)^3$.
	\end{lemma}
	\begin{scproof}
		At most $\Delta$ elements of $R^\shortleftarrow(y, \phi)$ are neighbors of $y$, and for each \ep{unordered} pair of colors $\alpha$, $\beta$, there are at most two vertices $z$ that are $(\phi, \alpha\beta)$-related to $y$ and satisfy $\deg(z, \phi, \alpha\beta) < 2$. Each such vertex $z$ has at most $\Delta$ neighbors, so
		\[
		|R^\shortleftarrow(y, \phi)| \,\leq\, 1 + \Delta + 2 \Delta {\Delta + 1 \choose 2} \,=\, \Delta^3 + \Delta^2 + \Delta + 1 \,\leq\, (\Delta+1)^3. \qedhere
		\]
	\end{scproof}

	Given a chain $C$, we write $R^\shortleftarrow(C, \phi) \defeq \bigcup_{y \in V(C)} R^\shortleftarrow(y, \phi)$.
	
	\begin{lemma}\label{lemma:change}
		Let $\phi$ be a proper partial coloring and let $C$ be a $\phi$-shiftable chain. Set $\psi \defeq \Shift(\phi, C)$. If $x$ is a vertex that is not in $R^\shortleftarrow(C, \phi) \cap R^\shortleftarrow(C, \psi)$, then $R^\shortrightarrow(x, \phi) = R^\shortrightarrow(x, \psi)$.
	\end{lemma}
	\begin{scproof}
		It is enough to show that if $x \not \in R^\shortleftarrow(C, \phi)$, then $R^\shortrightarrow(x, \phi) = R^\shortrightarrow(x, \psi)$; the conclusion for $x \not \in R^\shortleftarrow(C, \psi)$ then follows by exchanging the roles of $\phi$ and $\psi$ and replacing $C$ by its reverse $C^\ast$. So, suppose that $x \not \in R^\shortleftarrow(C, \phi)$. We first prove that $R^\shortrightarrow(x, \phi) \subseteq R^\shortrightarrow(x, \psi)$. Take any vertex $y \in R^\shortrightarrow(x, \phi)$. If $y=x$ or $y$ is adjacent to $x$, then $y \in R^\shortrightarrow(x, \psi)$ by definition, so suppose instead that there exist colors $\alpha$, $\beta$ and a neighbor $z$ of $x$ such that $\deg(z, \phi, \alpha\beta) < 2$ and $z$ and $y$ are $(\phi, \alpha\beta)$-related. In other words, $G(z, \phi, \alpha\beta)$ is a path starting at $z$ and containing $y$. Since $x \not \in R^\shortleftarrow(C, \phi)$, no vertex on this path is incident to an edge in $C$. This implies that $G(z, \psi, \alpha\beta) = G(z, \phi, \alpha\beta)$, so $G(z, \psi, \alpha\beta)$ is still a path starting at $z$ and passing through $y$, and hence $y \in R^\shortrightarrow(x, \psi)$, as desired.
		
		Now we show $R^\shortrightarrow(x, \psi) \subseteq R^\shortrightarrow(x, \phi)$. Take a vertex $y \in R^\shortrightarrow(x, \psi)$. Again, if $y = x$ or $y$ is adjacent to $x$, then we are done, so assume that there exist colors $\alpha$, $\beta$ and a neighbor $z$ of $x$ such that $\deg(z, \psi, \alpha\beta) < 2$ and $z$ and $y$ are $(\psi, \alpha\beta)$-related. Note that $\deg(z, \phi, \alpha\beta) < 2$ as well; otherwise, some edge incident to $z$ would be colored differently in $\phi$ and $\psi$, which would imply that $z \in V(C)$ and thus $x \in R^\shortleftarrow(C, \phi)$. Therefore, the graph $G(z, \phi, \alpha\beta)$ is a path starting at $z$. As before, we observe that no vertex on this path can be incident to an edge in $C$, since $x \not \in R^\shortleftarrow(C, \phi)$. But this means that $G(z, \psi, \alpha\beta) = G(z, \phi, \alpha\beta)$, so if $z$ and $y$ are $(\psi, \alpha\beta)$-related, then they are also $(\phi, \alpha\beta)$-related, as desired.
	\end{scproof}
	
	\subsection{The algorithm}\label{subsec:alg}
	
	For the remainder of \S\ref{sec:const}, we fix the following data:
	\begin{itemize}
		\item a proper partial coloring $\phi$;
		
		\item an uncolored edge $e = xy \in E \setminus \dom(\phi)$.
	\end{itemize}
	We also fix a numerical parameter $\ell \in \N^+$. \ep{We will eventually make $\ell$ be a value of order $\mathrm{poly}(\Delta)\log n$, but for now this does not matter.}
	We shall describe a randomized procedure, called the \emphd{Multi-Step Vizing Algorithm}, that attempts to build a $\phi$-happy chain $C$ with starting edge $e$. The chain $C$ will have the form $C = F_0 + P_0 + F_1 + P_1 + \cdots$, where each $F_i$ is a fan and each $P_i$ is a path-chain. The construction proceeds in steps, indexed by natural numbers \ep{the total number of steps may be finite or infinite}. On Step $i$, we build the fan $F_i$ and the path $P_i$.
	
	\subsubsection*{Initial set-up}
	
	Initially, we set $\phi_0 \defeq \phi$, $e_0 \defeq e$, $x_0 \defeq x$, and $y_0 \defeq y$. 
	
	\subsubsection*{Input for Step $i$}
	
	At the start of Step $i$, we have already constructed a number of interrelated structures. However, in order to describe how Step $i$ operates, only the following objects will be needed: 
	\begin{enumerate}[label=\ep{\normalfont{}In\arabic*}]
		\item a sequence of proper partial colorings $\phi_0$, $\psi_0$, \ldots, $\phi_{i-1}$, $\psi_{i-1}$, $\phi_i$;
		
		\item an edge $e_i = x_i y_i$ with $\phi_i(e_i) = \blank$;
		
		\item\label{item:colors} if $i > 0$, a pair of distinct colors $\alpha_{i-1}$ and $\beta_{i-1}$ such that $\deg(x_i, \phi_i, \alpha_{i-1}\beta_{i-1}) = 1$;
		
		\item a sequence $F_0$, $P_0$, \ldots, $F_{i-1}$, $P_{i-1}$, where each $F_j$ is a fan and each $P_j$ is a path-chain.
	\end{enumerate}
	
	\subsubsection*{Step $i$}
	
	We start by constructing a fan $F_i$. To this end, we consider two cases.
	
	\setcounter{case}{0}
	
	\begin{case}\label{case:happy_fan}
		There is a $\phi_i$-happy fan $F$ with $\pivot(F) = x_i$ and $\start(F) = e_i$.
	\end{case}
	
	\noindent Then we make $F_i$ be any such fan and stop the construction. In this case, the algorithm \emphd{succeeds}.
	
	\begin{case}\label{case:path}
		There is no $\phi_i$-happy fan $F$ with $\pivot(F) = x_i$ and $\start(F) = e_i$.
	\end{case}
	
	\noindent If $i = 0$, then we use Lemma~\ref{lemma:simple_fan} to obtain a $\phi_0$-shiftable fan $F_0$ with pivot $x_0$ and $\start(F_0) = e_0$ and a pair of distinct colors $\alpha_0$, $\beta_0$ such that $F_0$ is $(\phi_0, \alpha_0 \beta_0)$-successful. If, on the other hand, $i > 0$, then we apply Lemma~\ref{lemma:fancy_fan} instead with $\alpha_{i-1}$, $\beta_{i-1}$ in place of $\alpha$, $\beta$ \ep{this is possible by \ref{item:colors}} and obtain a $\phi_i$-shiftable fan $F_i$ with pivot $x_i$ and $\start(F_i) = e_i$ and a pair of distinct colors $\alpha_i$, $\beta_i$ such that:
	\begin{enumerate}[label=\ep{\normalfont{}F\arabic*}]
		\item\label{item:diff} either $\set{\alpha_i, \beta_i} \cap \set{\alpha_{i-1}, \beta_{i-1}} = \0$ and $F_i$ is $(\phi_i, \alpha_i \beta_i)$-successful,
		
		\item\label{item:same} or the following statements are valid:
		\begin{enumerate}[label=\ep{\normalfont{}F2\alph*}]
			\item $\set{\alpha_i, \beta_i} = \set{\alpha_{i-1}, \beta_{i-1}}$;
			\item $F_i$ is $(\phi_i, \alpha_i \beta_i)$-hopeful;
			\item\label{item:not_y} $\vend(F_i) \neq y_i$; and
			\item\label{item:ab_not_change} no edge in $F_i$ is colored $\alpha_{i}$ or $\beta_{i}$ by $\phi_i$.
		\end{enumerate}
	\end{enumerate}
	Note that regardless of whether \ref{item:diff} or \ref{item:same} holds, $F_i$ is $(\phi_i, \alpha_i\beta_i)$-hopeful. Let
	\[
	z_i \,\defeq\, \vend(F_i) \quad \text{and} \quad \psi_i \,\defeq\, \Shift(\phi_i, F_i).
	\]
	By Fact~\ref{fact:prom}, the edge $x_i z_i$ is $(\psi_i, \alpha_i\beta_i)$-hopeful, so we can define
	\[
	Q_i \,\defeq\, P(x_i, z_i, \psi_i, \alpha_i\beta_i).
	\]
	Set $u_i(0) \defeq x_i$, $u_i(1) \defeq z_i$ and write $Q_i = (u_i(0) u_i(1), u_i(1) u_i(2), u_i(2) u_i(3), \ldots)$.
	Next we proceed to define a path $P_i$, which will be an initial segment of $Q_i$. There are again a few cases to consider.
	
	\begin{subcase}\label{subcase:happy_path}
		The fan $F_i$ is $(\phi_i, \alpha_i\beta_i)$-successful and $\length(Q_i) \leq \ell$.
	\end{subcase}
	
	\noindent Then we let $P_i \defeq Q_i$ and stop the construction. In this case, the algorithm \emphd{succeeds}.
	
	\begin{subcase}\label{subcase:main}
		The fan $F_i$ is $(\phi_i, \alpha_i\beta_i)$-successful and $\length(Q_i) \geq \ell + 1$.
	\end{subcase}
	
	\noindent Note that since $\length(Q_i) \geq \ell+1$, the vertex $u_i(\ell+1)$ is well-defined. Pick an integer $1 \leq \ell_i \leq \ell$ uniformly at random and set
	\[
	x_{i+1} \,\defeq\,  u_i(\ell_i), \quad y_{i+1} \,\defeq\, u_i(\ell_i + 1), \quad \text{and} \quad P_i \,\defeq\, \rest{Q_i}{(\ell_i + 1)}.
	\]
	By construction, $P_i$ is a path-chain of length $\ell_i + 1 \leq \ell+1$ and we have
	\[
	\vend(P_i) \,=\, y_{i+1} \quad \text{and} \quad e_{i+1} \,\defeq\, x_{i+1} y_{i+1} \,=\, \eend(P_i).
	\]
	There are now two more subcases.
	
	\begin{subsubcase}\label{subsubcase:fail}
		At least one of the following holds:
		\begin{itemize}
			\item $x_{i+1} \in R^\shortleftarrow(F_j, \psi_j)$ for some $j \leq i$,
			
			\item or $x_{i+1} \in R^\shortleftarrow(P_j, \psi_j)$ for some $j \leq i - 1$.
		\end{itemize}
	\end{subsubcase}
	
	\noindent In this case we stop the construction and declare that the algorithm \emphd{fails}.
	
	\begin{subsubcase}\label{subsubcase:move_on}
		The following statements are valid:
		\begin{itemize}
			\item $x_{i+1} \not \in R^\shortleftarrow(F_j, \psi_j)$ for all $j \leq i$; and
			
			\item $x_{i+1} \not\in R^\shortleftarrow(P_j, \psi_j)$ for all $j \leq i - 1$.
		\end{itemize}
	\end{subsubcase}
	
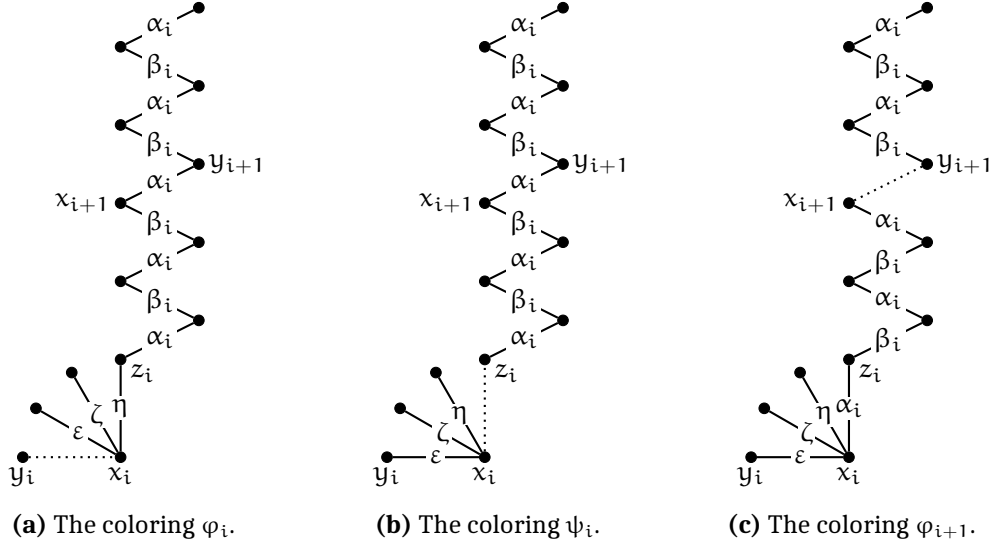
\begin{figure}[t]
		\centering
		\begin{subfigure}[t]{.2\textwidth}
			\centering
			\begin{tikzpicture}[scale=1.3]
			\node[circle,fill=black,draw,inner sep=0pt,minimum size=4pt] (a) at (0,0) {};
			\node[circle,fill=black,draw,inner sep=0pt,minimum size=4pt] (b) at (1,0) {};
			\path (b) ++(150:1) node[circle,fill=black,draw,inner sep=0pt,minimum size=4pt] (c) {}; 
			\path (b) ++(120:1) node[circle,fill=black,draw,inner sep=0pt,minimum size=4pt] (d) {}; 
			\node[circle,fill=black,draw,inner sep=0pt,minimum size=4pt] (e) at (1,1) {};
			\node[circle,fill=black,draw,inner sep=0pt,minimum size=4pt] (f) at (1.8,1.4) {};
			\node[circle,fill=black,draw,inner sep=0pt,minimum size=4pt] (g) at (1,1.8) {};
			\node[circle,fill=black,draw,inner sep=0pt,minimum size=4pt] (h) at (1.8,2.2) {};
			\node[circle,fill=black,draw,inner sep=0pt,minimum size=4pt] (i) at (1,2.6) {};
			\node[circle,fill=black,draw,inner sep=0pt,minimum size=4pt] (j) at (1.8,3) {};
			\node[circle,fill=black,draw,inner sep=0pt,minimum size=4pt] (k) at (1,3.4) {};
			\node[circle,fill=black,draw,inner sep=0pt,minimum size=4pt] (l) at (1.8,3.8) {};
			\node[circle,fill=black,draw,inner sep=0pt,minimum size=4pt] (m) at (1,4.2) {};
			\node[circle,fill=black,draw,inner sep=0pt,minimum size=4pt] (n) at (1.8,4.6) {};
			
			\draw[ thick,dotted] (a) -- (b);
			\draw[ thick] (b) to node[midway,inner sep=1pt,outer sep=1pt,minimum size=4pt,fill=white] {$\epsilon$} (c) (b) to node[midway,inner sep=0.5pt,outer sep=0.5pt,minimum size=4pt,fill=white] {$\zeta$} (d) (b) to node[midway,inner sep=1pt,outer sep=1pt,minimum size=4pt,fill=white] {$\eta$} (e) to node[midway,inner sep=1pt,outer sep=1pt,minimum size=4pt,fill=white] {$\alpha_i$} (f) to node[midway,inner sep=1pt,outer sep=1pt,minimum size=4pt,fill=white] {$\beta_i$} (g) to node[midway,inner sep=1pt,outer sep=1pt,minimum size=4pt,fill=white] {$\alpha_i$} (h) to node[midway,inner sep=1pt,outer sep=1pt,minimum size=4pt,fill=white] {$\beta_i$} (i) to node[midway,inner sep=1pt,outer sep=1pt,minimum size=4pt,fill=white] {$\alpha_i$} (j) to node[midway,inner sep=1pt,outer sep=1pt,minimum size=4pt,fill=white] {$\beta_i$} (k) to node[midway,inner sep=1pt,outer sep=1pt,minimum size=4pt,fill=white] {$\alpha_i$} (l) to node[midway,inner sep=1pt,outer sep=1pt,minimum size=4pt,fill=white] {$\beta_i$} (m) to node[midway,inner sep=1pt,outer sep=1pt,minimum size=4pt,fill=white] {$\alpha_i$} (n);
			
			
			
			
			\node[anchor=148.28] at (e) {$z_i$};
			\node[anchor=north] at (a) {$y_i$};
			\node[anchor=north] at (b) {$x_i$};
			\node[anchor=east] at (i) {$x_{i+1}$};
			\node[anchor=west] at (j) {$y_{i+1}$};
			\end{tikzpicture}
			\caption{The coloring $\phi_i$.}
		\end{subfigure}%
		\qquad\qquad%
		\begin{subfigure}[t]{.2\textwidth}
			\centering
			\begin{tikzpicture}[scale=1.3]
			\node[circle,fill=black,draw,inner sep=0pt,minimum size=4pt] (a) at (0,0) {};
			\node[circle,fill=black,draw,inner sep=0pt,minimum size=4pt] (b) at (1,0) {};
			\path (b) ++(150:1) node[circle,fill=black,draw,inner sep=0pt,minimum size=4pt] (c) {}; 
			\path (b) ++(120:1) node[circle,fill=black,draw,inner sep=0pt,minimum size=4pt] (d) {}; 
			\node[circle,fill=black,draw,inner sep=0pt,minimum size=4pt] (e) at (1,1) {};
			\node[circle,fill=black,draw,inner sep=0pt,minimum size=4pt] (f) at (1.8,1.4) {};
			\node[circle,fill=black,draw,inner sep=0pt,minimum size=4pt] (g) at (1,1.8) {};
			\node[circle,fill=black,draw,inner sep=0pt,minimum size=4pt] (h) at (1.8,2.2) {};
			\node[circle,fill=black,draw,inner sep=0pt,minimum size=4pt] (i) at (1,2.6) {};
			\node[circle,fill=black,draw,inner sep=0pt,minimum size=4pt] (j) at (1.8,3) {};
			\node[circle,fill=black,draw,inner sep=0pt,minimum size=4pt] (k) at (1,3.4) {};
			\node[circle,fill=black,draw,inner sep=0pt,minimum size=4pt] (l) at (1.8,3.8) {};
			\node[circle,fill=black,draw,inner sep=0pt,minimum size=4pt] (m) at (1,4.2) {};
			\node[circle,fill=black,draw,inner sep=0pt,minimum size=4pt] (n) at (1.8,4.6) {};
			
			\draw[ thick] (a) to node[midway,inner sep=1pt,outer sep=1pt,minimum size=4pt,fill=white] {$\epsilon$} (b);
			\draw[ thick,dotted] (b) -- (e);
			\draw[ thick] (b) to node[midway,inner sep=1pt,outer sep=1pt,minimum size=4pt,fill=white] {$\zeta$} (c) (b) to node[midway,inner sep=0.5pt,outer sep=0.5pt,minimum size=4pt,fill=white] {$\eta$} (d)  (e) to node[midway,inner sep=1pt,outer sep=1pt,minimum size=4pt,fill=white] {$\alpha_i$} (f) to node[midway,inner sep=1pt,outer sep=1pt,minimum size=4pt,fill=white] {$\beta_i$} (g) to node[midway,inner sep=1pt,outer sep=1pt,minimum size=4pt,fill=white] {$\alpha_i$} (h) to node[midway,inner sep=1pt,outer sep=1pt,minimum size=4pt,fill=white] {$\beta_i$} (i) to node[midway,inner sep=1pt,outer sep=1pt,minimum size=4pt,fill=white] {$\alpha_i$} (j) to node[midway,inner sep=1pt,outer sep=1pt,minimum size=4pt,fill=white] {$\beta_i$} (k) to node[midway,inner sep=1pt,outer sep=1pt,minimum size=4pt,fill=white] {$\alpha_i$} (l) to node[midway,inner sep=1pt,outer sep=1pt,minimum size=4pt,fill=white] {$\beta_i$} (m) to node[midway,inner sep=1pt,outer sep=1pt,minimum size=4pt,fill=white] {$\alpha_i$} (n);
			
			
			
			
			\node[anchor=148.28] at (e) {$z_i$};
			\node[anchor=north] at (a) {$y_i$};
			\node[anchor=north] at (b) {$x_i$};
			\node[anchor=east] at (i) {$x_{i+1}$};
			\node[anchor=west] at (j) {$y_{i+1}$};
			\end{tikzpicture}
			\caption{The coloring $\psi_{i}$.}\label{fig:just_shifted:2}
		\end{subfigure}%
		\qquad\qquad%
		\begin{subfigure}[t]{.2\textwidth}
			\centering
			\begin{tikzpicture}[scale=1.3]
			\node[circle,fill=black,draw,inner sep=0pt,minimum size=4pt] (a) at (0,0) {};
			\node[circle,fill=black,draw,inner sep=0pt,minimum size=4pt] (b) at (1,0) {};
			\path (b) ++(150:1) node[circle,fill=black,draw,inner sep=0pt,minimum size=4pt] (c) {}; 
			\path (b) ++(120:1) node[circle,fill=black,draw,inner sep=0pt,minimum size=4pt] (d) {}; 
			\node[circle,fill=black,draw,inner sep=0pt,minimum size=4pt] (e) at (1,1) {};
			\node[circle,fill=black,draw,inner sep=0pt,minimum size=4pt] (f) at (1.8,1.4) {};
			\node[circle,fill=black,draw,inner sep=0pt,minimum size=4pt] (g) at (1,1.8) {};
			\node[circle,fill=black,draw,inner sep=0pt,minimum size=4pt] (h) at (1.8,2.2) {};
			\node[circle,fill=black,draw,inner sep=0pt,minimum size=4pt] (i) at (1,2.6) {};
			\node[circle,fill=black,draw,inner sep=0pt,minimum size=4pt] (j) at (1.8,3) {};
			\node[circle,fill=black,draw,inner sep=0pt,minimum size=4pt] (k) at (1,3.4) {};
			\node[circle,fill=black,draw,inner sep=0pt,minimum size=4pt] (l) at (1.8,3.8) {};
			\node[circle,fill=black,draw,inner sep=0pt,minimum size=4pt] (m) at (1,4.2) {};
			\node[circle,fill=black,draw,inner sep=0pt,minimum size=4pt] (n) at (1.8,4.6) {};
			
			\draw[ thick] (a) to node[midway,inner sep=1pt,outer sep=1pt,minimum size=4pt,fill=white] {$\epsilon$} (b);
			\draw[dotted, thick] (i) -- (j);
			\draw[ thick] (b) to node[midway,inner sep=0.5pt,outer sep=0.5pt,minimum size=4pt,fill=white] {$\zeta$} (c) (b) to node[midway,inner sep=1pt,outer sep=1pt,minimum size=4pt,fill=white] {$\eta$} (d) (b) to node[midway,inner sep=1pt,outer sep=1pt,minimum size=4pt,fill=white] {$\alpha_i$} (e) to node[midway,inner sep=1pt,outer sep=1pt,minimum size=4pt,fill=white] {$\beta_i$} (f) to node[midway,inner sep=1pt,outer sep=1pt,minimum size=4pt,fill=white] {$\alpha_i$} (g) to node[midway,inner sep=1pt,outer sep=1pt,minimum size=4pt,fill=white] {$\beta_i$} (h) to node[midway,inner sep=1pt,outer sep=1pt,minimum size=4pt,fill=white] {$\alpha_i$} (i)  (j) to node[midway,inner sep=1pt,outer sep=1pt,minimum size=4pt,fill=white] {$\beta_i$} (k) to node[midway,inner sep=1pt,outer sep=1pt,minimum size=4pt,fill=white] {$\alpha_i$} (l) to node[midway,inner sep=1pt,outer sep=1pt,minimum size=4pt,fill=white] {$\beta_i$} (m) to node[midway,inner sep=1pt,outer sep=1pt,minimum size=4pt,fill=white] {$\alpha_i$} (n);
			
			
			
			\node[anchor=148.28] at (e) {$z_i$};
			\node[anchor=north] at (a) {$y_i$};
			\node[anchor=north] at (b) {$x_i$};
			\node[anchor=east] at (i) {$x_{i+1}$};
			\node[anchor=west] at (j) {$y_{i+1}$};
			\end{tikzpicture}
			\caption{The coloring $\phi_{i+1}$.}\label{fig:just_shifted:3}
		\end{subfigure}%
		\caption{The colorings $\phi_i$, $\psi_i$, and $\phi_{i+1}$. In this case, $\ell_i = 5$.}\label{fig:just_shifted}
	\end{figure}

	\noindent Then we let $\phi_{i+1} \defeq \Shift(\psi_i, P_i)$. Note that the unique edge of color $\alpha_i$ or $\beta_i$ in $\phi_{i+1}$ incident to $x_{i+1}$ is $u_i(\ell_i - 1) x_{i+1}$, so $\deg(x_{i+1}, \phi_{i+1}, \alpha_i \beta_i) = 1$ \ep{see Fig.~\ref{fig:just_shifted:3}}. Thus, \ref{item:colors} holds with $i+1$ in place of $i$, and hence we may continue to Step $i+1$. This finishes Subcase~\ref{subcase:main}.
	
	\begin{samepage}
		\begin{subcase}\label{subcase:impossible}
			The fan $F_i$ is not $(\phi_i, \alpha_i\beta_i)$-successful.
		\end{subcase}
		
		\noindent We will show that this situation cannot occur \ep{see Lemma~\ref{lemma:impossible}}. For now let us say that in this case the procedure stops and the algorithm \emphd{fails}.
	\end{samepage}
	
	\subsubsection*{Output of Step $i$}
	
	Here is a brief r\'esum\'e of the possible outcomes for Step $i$:
	
	\begin{enumerate}[label=\ep{\normalfont{}Out\arabic*}]
		\item In Case~\ref{case:happy_fan}, the algorithm successfully  outputs a $\phi$-happy chain $F_0 + P_0 + \cdots + P_{i-1} + F_i$.
		
		\item In Subcase~\ref{subcase:happy_path}, the algorithm successfully outputs a $\phi$-happy chain $F_0 + P_0 + \cdots + F_i + P_i$.
		
		\item In Subsubcase~\ref{subsubcase:fail}, the algorithm fails.
		
		\item In Subsubcase~\ref{subsubcase:move_on}, the algorithm moves on to Step $i+1$.
		
		\item In Subcase~\ref{subcase:impossible}, the algorithm fails, but we will prove in Lemma~\ref{lemma:impossible} that this cannot happen.
	\end{enumerate}
	Note that the only case in which the algorithm moves on to the next step is Subsubcase~\ref{subsubcase:move_on}.
	
	\subsection{The backtracking lemma}
	
	In this subsection we prove a lemma that plays a central role in our analysis of the Multi-Step Vizing Algorithm. In essence, it formalizes the idea of tracing back alternating paths that we sketched in \S\ref{sec:overview}.
	
	\begin{lemma}[\textls{Backtracking}]\label{lemma:back_orbit}
		Suppose that on Step $i$ of the Multi-Step Vizing Algorithm, Case~\ref{case:path} occurs. Fix a vertex $u \in V(Q_i)$. Then $x_i \in R^\shortleftarrow(u, \psi_j) \cap R^\shortleftarrow(u, \phi_j)$ for all $j \leq i$, and, in particular, $x_i \in R^\shortleftarrow(u, \phi)$.
	\end{lemma}
	\begin{scproof}\stepcounter{ForClaims} \renewcommand{\theForClaims}{\ref{lemma:back_orbit}}
		Let $\alpha \defeq \alpha_i$ and $\beta \defeq \beta_i$. If $u = x_i$, then we are done, so assume that $u \neq x_i$. Then $x_i \in R^\shortleftarrow(u, \psi_i)$ since $z_i$ is a neighbor of $x_i$ with $\deg(z_i, \psi_i, \alpha\beta) < 2$ and $u$ is $(\psi_i, \alpha\beta)$-related to $z_i$. Next we show that
		\begin{equation}\label{eq:phi_i}
		x_i \in R^\shortleftarrow(u, \phi_i, \alpha\beta).
		\end{equation}
		Since the fan $F_i$ is $(\phi_i, \alpha\beta)$-hopeful, we have $\deg(z_i, \phi_i, \alpha\beta) < 2$. Suppose, towards a contradiction, that $u$ and $z_i$ are not $(\phi_i, \alpha\beta)$-related. By construction, $u$ and $z_i$ are joined by a path in $G(\psi_i, \alpha\beta)$. If $u$ and $z_i$ are not $(\phi_i, \alpha\beta)$-related, then at least one edge on this path must have a color different from $\alpha$ and $\beta$ in $\phi_i$, which means that it must belong to $F_i$. But every edge in $F_i$ contains $x_i$, so $x_i$ is an internal vertex on some path in $G(\psi_i, \alpha\beta)$, which is impossible as $\deg(x_i, \psi_i, \alpha\beta) =  \deg(x_i, \phi_i, \alpha\beta) < 2$.
		
		If $i =0$, then we are done, so assume that $i > 0$. On Step $i-1$, the algorithm must have entered Subsubcase~\ref{subsubcase:move_on} \ep{otherwise it would have stopped without reaching Step $i$}, which means that the following statements hold:
		\begin{itemize}
			\item $x_i \not \in R^\shortleftarrow(F_j, \psi_j)$ for all $j \leq i-1$; and
			
			\item $x_i \not\in R^\shortleftarrow(P_j, \psi_j)$ for all $j \leq i - 2$.
		\end{itemize}
		By iteratively applying Lemma~\ref{lemma:change}, we conclude that
		\[
		R^\shortrightarrow(x_i, \phi_0) \,=\, R^\shortrightarrow(x_i, \psi_0) \,=\, R^\shortrightarrow(x_i, \phi_1) \,=\, R^\shortrightarrow(x_i, \psi_1) \,=\, \cdots \,=\, R^\shortrightarrow(x_i, \psi_{i-1}).
		\]
		Thus, it remains to argue that $x_i \in R^\shortleftarrow(u, \psi_{i-1})$. We will in fact show that
		\[
		x_i \in R^\shortleftarrow(u, \psi_{i-1}, \alpha\beta).
		\]
		There are two possibilities to consider. If on Step $i$, situation \ref{item:diff} occurs, then $\set{\alpha, \beta} \cap \set{\alpha_{i-1}, \beta_{i-1}} = \0$. The only colors that are involved in the change from $\psi_{i-1}$ to $\phi_i$ are $\alpha_{i-1}$ and $\beta_{i-1}$, so in this case $G(\psi_{i-1}, \alpha\beta) = G(\phi_i, \alpha\beta)$. In particular, \eqref{eq:phi_i} yields $x_i \in R^\shortleftarrow(u, \psi_{i-1}, \alpha\beta)$, as desired. Now assume that on Step $i$, situation \ref{item:same} occurs and, in particular, $\set{\alpha_{i-1}, \beta_{i-1}} = \set{\alpha, \beta}$.
		
		\begin{claim*}
			Under the above assumptions, $\deg(z_i, \psi_{i-1}, \alpha\beta) = \deg(z_i, \phi_i, \alpha\beta) < 2$.
		\end{claim*}
		\begin{claimproof}
			When the path $P_{i-1}$ is shifted, the only vertices whose degrees in the $\alpha\beta$-colored graph change are the first and the last two vertices of $P_{i-1}$, i.e., $x_{i-1}$, $z_{i-1}$, $x_i$, and $y_i$. Since $x_i$ and $z_i$ are adjacent, $z_i \neq x_i$. Also, $z_i \neq y_i$ by \ref{item:not_y}. The only remaining options are $z_i = x_{i-1}$ and $z_i = z_{i-1}$, but both of them imply $x_i \in R^\shortleftarrow(F_{i-1}, \psi_{i-1})$, which is a contradiction.
		\end{claimproof}
		
		It remains to check that $u$ and $z_i$ are $(\psi_{i-1}, \alpha\beta)$-related. Suppose, towards a contradiction, that they are not. This means that the subgraphs $G(z_i, \psi_{i-1}, \alpha\beta)$ and $G(z_i, \phi_i, \alpha \beta)$ are distinct, and hence at least one vertex in $G(z_i, \psi_{i-1}, \alpha \beta)$ is incident to an edge that is colored differently in $\psi_{i-1}$ and in $\phi_i$. All such vertices, except $x_{i-1}$, are $(\psi_{i-1}, \alpha \beta)$-related to $z_{i-1}$ \ep{see Fig.~\ref{fig:just_shifted:2} for an illustration}. Thus, $z_i$ is $(\psi_{i-1}, \alpha \beta)$-related either to $x_{i-1}$ or to $z_{i-1}$. In any case, $x_i \in R^\shortleftarrow(F_{i-1}, \psi_{i-1})$; a contradiction.
	\end{scproof}

	Next we apply Lemma~\ref{lemma:back_orbit} to show that Subcase~\ref{subcase:impossible} cannot occur:
	
	\begin{lemma}\label{lemma:impossible}
		On Step $i$ of the Multi-Step Vizing Algorithm, Subcase~\ref{subcase:impossible} cannot occur.
	\end{lemma}
	\begin{scproof}
		Let $\alpha \defeq \alpha_i$ and $\beta \defeq \beta_i$ and suppose, towards a contradiction, that the fan $F_i$ is not $(\phi_i, \alpha\beta)$-successful. This implies that $i > 0$ and \ref{item:same} holds; in particular, $\set{\alpha_{i-1}, \beta_{i-1}} = \set{\alpha, \beta}$. Since $F_i$ is $(\phi_i, \alpha\beta)$-hopeful but not $(\phi_i, \alpha\beta)$-successful, the vertices $x_i$ and $z_i$ are $(\psi_i, \alpha\beta)$-related. From \ref{item:ab_not_change}, we conclude that $G(\psi_i, \alpha\beta) = G(\phi_i, \alpha\beta)$, so $x_i$ and $z_i$ are also $(\phi_i, \alpha\beta)$-related. Since $\set{\alpha_{i-1}, \beta_{i-1}} = \set{\alpha, \beta}$, $x_i$ is $(\phi_i, \alpha\beta)$-related to $x_{i-1}$ \ep{see Fig.~\ref{fig:just_shifted:3} for an illustration}. Therefore, $z_i$ and $x_{i-1}$ are $(\phi_i, \alpha\beta)$-related, and, by \ref{item:ab_not_change} again, $z_i$ and $x_{i-1}$ are $(\psi_i, \alpha\beta)$-related. This implies that $x_{i-1} \in V(Q_i)$, so, by Lemma~\ref{lemma:back_orbit}, $x_i \in R^\shortleftarrow(x_{i-1}, \psi_{i-1})$ and thus $x_i \in R^\shortleftarrow(F_{i-1}, \psi_{i-1})$. But if this were the case, then the algorithm would have entered Subsubcase~\ref{subsubcase:fail} and stopped on Step $i-1$; a contradiction.
	\end{scproof}

	\section{Probabilistic analysis of the Multi-Step Vizing Algorithm}\label{sec:analysis}
	
	\subsection{The algorithm is likely to succeed}
	
	We use the set-up and the notation of \S\ref{subsec:alg}. Our goal is to obtain a lower bound on the probability that the Multi-Step Vizing Algorithm terminates successfully within a limited number of steps.
	
	\begin{lemma}\label{lemma:lucky}
		Fix $T \in \N^+$ and set $\lambda \defeq \ell /(\Delta + 1)^3$. Assuming that $\lambda > 1$, the probability that the Multi-Step Vizing Algorithm successfully terminates on Step $i$ for some $i < T$ is at least
		\[
		1 \,-\, \frac{n}{\lambda^T} \,-\, \frac{3T(\Delta+1)^3}{\lambda - 1}.
		\]
	\end{lemma}
	\begin{scproof}\stepcounter{ForClaims} \renewcommand{\theForClaims}{\ref{lemma:lucky}}
		First we bound the probability that the algorithm enters Step $T$.
		Fix a vertex $u \in V$ and consider the probability $\P\left[x_T = u\right]$. Notice that if $x_T = u$, then, by Lemma~\ref{lemma:back_orbit},
		\begin{equation}\label{eq:backtracking}
		x_T = u, \quad x_{T-1} \in R^\shortleftarrow(x_T, \phi), \quad x_{T-2} \in R^\shortleftarrow(x_{T-1}, \phi), \quad \ldots, \quad \text{and} \quad x_1 \in R^\shortleftarrow(x_2, \phi).
		\end{equation}
		By Lemma~\ref{lemma:indegree}, the number of sequences $(x_1, \ldots, x_T)$ satisfying \eqref{eq:backtracking} is at most $(\Delta + 1)^{3(T-1)}$. Since each of the vertices $x_1$, \ldots, $x_T$ is chosen uniformly at random from a set of $\ell$ candidates, we conclude that
		\[
		\P\left[x_T = u\right] \,\leq\, \frac{(\Delta+1)^{3(T-1)}}{\ell^T} \,\leq\, \frac{1}{\lambda^T}.
		\]
		Since the vertex $x_T$ must be defined for the algorithm to enter Step $T$, we have
		\[
		\P\left[\text{the algorithm enters Step $T$}\right] \,\leq\, \sum_{u \in V} \P\left[x_T = u\right] \,\leq\, \frac{n}{\lambda^T}.
		\]
		
		Next we bound the probability that the algorithm fails on Step $i$ for some $i < T$. Thanks to Lemma~\ref{lemma:impossible}, the only situation in which the algorithm fails on Step $i$ is Subsubcase~\ref{subsubcase:fail}, i.e.,
		\begin{itemize}
			\item $x_{i+1} \in R^\shortleftarrow(F_j, \psi_j)$ for some $j \leq i$,
			
			\item or $x_{i+1} \in R^\shortleftarrow(P_j, \psi_j)$ for some $j \leq i - 1$.
		\end{itemize}
		
		\begin{claim*}
			For all $j \leq i$, we have
			\[
			\P\left[x_{i+1} \in R^\shortleftarrow(F_j, \psi_j)\right] \,\leq\, \frac{\Delta+1}{\lambda^{i-j+1}} \quad \text{and} \quad \P\left[x_{i+1} \in R^\shortleftarrow(P_j, \psi_j)\right] \,\leq\, \frac{2(\Delta+1)^{3}}{\lambda^{i-j}}.
			\]
		\end{claim*}
		\begin{claimproof}
			We will show that the desired upper bounds hold even if we allow an adversary to specify the 
			outcome of the algorithm's execution until it is time to pick $\ell_j$ \ep{if the algorithm terminates earlier, or if it does not enter Subcase~\ref{subcase:main} on Step $j$, then the events $x_{i+1} \in R^\shortleftarrow(F_j, \psi_j)$ and $x_{i+1} \in R^\shortleftarrow(P_j, \psi_j)$ cannot happen, so we may ignore that possibility}. So, fix the outcome of the algorithm's execution up to the stage when $\ell_j$ is supposed to be chosen and let $\P^\ast$ denote conditional probability subject to this outcome. At this point, we have already constructed the fan $F_j$, the coloring $\psi_j$, and the path $Q_j = (u_j(0) u_j(1), u_j(1) u_j(2), \ldots)$ of length at least $\ell + 1$, but we have not yet picked the vertex $x_{j+1}$. Fix any $u \in V$ and consider the \ep{conditional} probability $\P^\ast\left[x_{i+1} = u\right]$. If $x_{i+1} = u$, then, by Lemma~\ref{lemma:back_orbit},
			\begin{equation}\label{eq:backtracking1}
			x_{i+1} = u, \quad x_i \in R^\shortleftarrow(x_{i+1}, \phi), \quad x_{i-1} \in R^\shortleftarrow(x_i, \phi), \quad \ldots, \quad \text{and} \quad x_{j+1} \in R^\shortleftarrow(x_{j+2}, \phi).
			\end{equation}
			By Lemma~\ref{lemma:indegree}, the number of sequences $(x_{j+1}, \ldots, x_{i+1})$ satisfying \eqref{eq:backtracking1} is at most $(\Delta + 1)^{3(i-j)}$. Since each of the vertices $x_{j+1}$, \ldots, $x_{i+1}$ is chosen uniformly at random from a set of $\ell$ candidates, we get
			\[
			\P^\ast\left[x_{i+1} = u\right] \,\leq\, \frac{(\Delta+1)^{3(i-j)}}{\ell^{i-j+1}} \,=\, \frac{1}{ \lambda^{i-j}\ell}.
			\]
			Set $U \defeq R^\shortleftarrow(F_j, \psi_j)$. By Lemma~\ref{lemma:indegree}, we have $|U| \leq (\Delta+1)^3 |V(F_j)| \leq (\Delta+1)^4$, so
			\[
			\P^\ast\left[x_{i+1} \in R^\shortleftarrow(F_j, \psi_j)\right] \,=\, \sum_{u \in U} \P^\ast\left[x_{i+1} = u\right] \,\leq\, \frac{(\Delta+1)^4}{\lambda^{i-j}\ell} \,=\, \frac{\Delta+1}{\lambda^{i-j+1}}.
			\]
			Regardless of the value $\ell_j$, the set $R^\shortleftarrow(P_j, \psi_j)$ is contained in $W \defeq \bigcup_{k = 0}^{\ell + 1} R^\shortleftarrow(u_j(k), \psi_j)$. By Lemma~\ref{lemma:indegree}, $|W| \leq (\Delta+1)^3 (\ell + 2) \leq 2(\Delta+1)^3 \ell$ \ep{where we are using that $\ell > (\Delta+1)^3 \geq 1$, so $\ell+2 \leq 2\ell$}, and hence
			\[
			\P^\ast\left[x_{i+1} \in R^\shortleftarrow(P_j, \psi_j)\right] \,\leq\, \sum_{u \in W} \P^\ast\left[x_{i+1} = u\right] \,\leq\, \frac{2(\Delta+1)^3}{\lambda^{i-j}}.\qedhere
			\]
		\end{claimproof}
		
		Using the above claim and the union bound, we obtain
		\begin{align*}
		\P[\text{the }&\text{algorithm fails on Step $i$ for some $i < T$}] \\
		&\leq\, \sum_{j=0}^{T-1} \sum_{i= j}^\infty \P\left[x_{i+1} \in R^\shortleftarrow(F_j, \psi_j)\right] \,+\, \sum_{j=0}^{T-2} \sum_{i=j+1}^\infty \P\left[x_{i+1} \in R^\shortleftarrow(P_j, \psi_j)\right] \\
		[\text{by the claim}] \qquad &\leq\, \sum_{j=0}^{T-1} \sum_{i= j}^\infty \frac{\Delta+1}{\lambda^{i-j+1}} \,+\, \sum_{j=0}^{T-2} \sum_{i=j+1}^\infty \frac{2(\Delta+1)^{3}}{\lambda^{i-j}} \\
		&=\, \sum_{j=0}^{T-1} \sum_{k=1}^\infty \frac{\Delta+1}{\lambda^{k}}  \,+\, \sum_{j=0}^{T-2} \sum_{k=1}^\infty \frac{2(\Delta+1)^{3}}{\lambda^{k}} \\
		[\text{since $\lambda > 1$}] \qquad &=\, \frac{T(\Delta+1)}{\lambda-1} \,+\, \frac{2(T-1)(\Delta+1)^3}{\lambda-1} \,\leq\, \frac{3T(\Delta+1)^3}{\lambda - 1}.
		\end{align*}
		Putting everything together, we conclude that
		\begin{align*}
		\P[\text{the }&\text{algorithm succeeds on Step $i$ for some $i < T$}] \\
		&=\,1 \,-\, \P\left[\text{the algorithm enters Step $T$}\right] \,-\, \P[\text{the algorithm fails on Step $i$ for some $i < T$}] \\
		&\leq\, 1 \,-\, \frac{n}{\lambda^T} \,-\, \frac{3T(\Delta+1)^3}{\lambda - 1}.\qedhere
		\end{align*}
	\end{scproof}
	
	Lemma~\ref{lemma:lucky} immediately yields the existence of ``small'' augmenting subgraphs:
	
	\begin{corl}[to Lemma~\ref{lemma:lucky}]
		Assuming $n$ is sufficiently large, the following statement holds: For every proper partial coloring $\phi$ and every uncolored edge $e = xy$, there is a connected augmenting subgraph $H$ whose only uncolored edge is $e$ with $|E(H)| \leq (\Delta+1)^6(\log n)^2$.
	\end{corl}
	\begin{scproof}
		By Fact~\ref{fact:happy}, it suffices to find a $\phi$-happy chain $C$ with $\start(C) = e$ and \[\length(C) \,\leq\, (\Delta+1)^6(\log n)^2.\] To this end, fix a small positive constant $\epsilon$ \ep{any $\epsilon < 1/3$ would work}, let $T \defeq \lfloor \epsilon \log n \rfloor$, and run the Multi-Step Vizing Algorithm with $e_0 = e$ and
		\[
		\ell \,\defeq\, 
		(\Delta+1)^6 \lfloor \log n \rfloor.
		\]
		Since $\lambda \defeq \ell/(\Delta+1)^3 = (\Delta+1)^3\lfloor \log n \rfloor > 1$, we can use Lemma~\ref{lemma:lucky} to conclude that the algorithm successfully terminates on Step $i$ for some $i < T$ with probability at least
		\[
		1 \,-\, \frac{n}{\lambda^T} \,-\, \frac{3T(\Delta+1)^3}{\lambda - 1} \,=\, 1 \,-\, n \left(\frac{1}{(\Delta+1)^3\lfloor \log n \rfloor}\right)^{\lfloor \epsilon \log n \rfloor}\,-\, 3\epsilon \,-\, o(1) \,=\, 1 - 3\epsilon - o(1) \,>\, 0.
		\]
		Whenever the algorithm terminates successfully on Step $i < T$, it produces a $\phi$-happy chain $C$ with $\start(C) = e$ of the form $C = F_0 + P_0 + \cdots + P_{i-1} + F_i$ or $C = F_0 + P_0 + \cdots + F_i + P_i$. Since $\length(F_i) \leq \Delta$ and $\length(P_i) \leq \ell + 1$ for all $i$, we obtain
		\[
		\length(C) \,\leq\, T(\ell + \Delta) \,\leq\, (\epsilon + o(1)) (\Delta+1)^6(\log n)^2,
		\]
		which yields the desired bound.
	\end{scproof}
	
	\subsection{Proof of Theorem~\ref{theo:ext}}
	
	For the reader's convenience, we recall the statement of Theorem~\ref{theo:ext}.
	
	\begin{theocopy}{theo:ext}
		There is $n_0 \in \N$ such that for all $n > n_0$, the following holds. Let $G$ be an $n$-vertex graph of maximum degree $\Delta$. Fix a partition $E(G) = X \sqcup U$ and let $\phi \colon X \to [\Delta + 1]$ be a proper partial edge-coloring. Then there exists a subset $W \subseteq U$ of size $|W| \geq |U|/((\Delta+1)^{10}(\log n)^2)$ such that it is possible to assign to each edge $e \in W$ a connected augmenting subgraph $H^e \subseteq G$ with the following properties:
		\begin{itemize}
			\item for each $e \in W$, 
			$E(H^e) \cap U = \set{e}$, that is, $e$ is the unique uncolored edge in $H^e$;
			\item the vertex sets of the graphs $H^e$, $e \in W$, are pairwise disjoint;
			\item for each $e \in W$, $|E(H^e)| \leq (\Delta+1)^6(\log n)^2$.
		\end{itemize}
	\end{theocopy}
	\begin{scproof}
		Throughout the proof, we will be assuming that $n$ is sufficiently large. Let the endpoints of each uncolored edge $e$ be $x^e$ and $y^e$ \ep{it does not matter which endpoint we call $x^e$ and which $y^e$}. Set the value of the parameter $\ell$ to
		\[
		\ell \,\defeq\, 
		(\Delta+1)^6 \lfloor\log n\rfloor
		\]
		and run the Multi-Step Vizing Algorithm with $e_0 = e$, $x_0 = x^e$, and $y_0 = y^e$ for every uncolored edge $e$, making the processes for different uncolored edges independent of each other. We use variables with ``$e$'' as a superscript, such as $F_i^e$ and $P_i^e$, to denote the structures obtained using the execution of the algorithm with $e_0 = e$. In particular, for each $e \in U$, we have a sequence of vertices $x^e = x_0^e$, $x_1^e$, $x_2^e$, \ldots.
		
		Let $\lambda \defeq \ell /(\Delta+1)^3 = (\Delta+1)^3\lfloor \log n \rfloor$. Fix a small positive constant $\epsilon$ \ep{any $\epsilon < 1/12$ would work} and set $T \defeq \lfloor \epsilon \log n \rfloor$. Say that an uncolored edge $e$ is \emphd{lucky} if the execution of the Multi-Step Vizing Algorithm with $e_0 = e$ successfully terminated on Step $i$ for some $i < T$. Let $L \subseteq U$ denote the set of all lucky edges. For $e \in L$, we let $C^e$ be the $\phi$-happy chain with $\start(C^e) = e$ produced by the algorithm with $e_0 = e$. In other words, if the algorithm terminated on Step $i$ due to Case~\ref{case:happy_fan}, then $C^e = F_0^e + P_0^e + \cdots + P_{i-1}^e + F_i^e$, and if it terminated due to Subcase~\ref{subcase:happy_path}, then $C^e = F_0^e + P_0^e + \cdots + F_i^e + P_i^e$. If $e \in L$, then we let $H^e \defeq H(C^e)$. By construction, $H^e$ is a connected $\phi$-augmenting subgraph whose only uncolored edge is $e$. Furthermore, since $e$ is lucky,
		\begin{equation}
		|E(H^e)|,\, |V(H^e)| \,\leq\, \length(C^e) + 1 \,\leq\, T(\ell+\Delta) + 1 \,\leq\, (\epsilon + o(1)) (\Delta+1)^6(\log n)^2.\label{eq:size_bound}
		\end{equation}
		
		Our goal now is to show that, with positive probability, there is a set $W \subseteq L$ of lucky edges of size at least $|W| \geq |U|/((\Delta+1)^{10}(\log n)^2)$ such that the graphs $H^e$, $e \in E$, are pairwise vertex-disjoint. To this end, define an auxiliary graph $\Gamma$ with $V(\Gamma) \defeq E$ by making two distinct edges $e$, $h$ adjacent in $\Gamma$ if and only if $e$ and $h$ are both lucky and $V(H^e) \cap V(H^h) \neq \0$.
		
		\begin{claim*}
			For each uncolored edge $e$, $\mathbb{E}\left[\deg_\Gamma(e)\right] \leq (\epsilon+o(1))(\Delta+1)^{10} (\log n)^2$.
		\end{claim*}
		\begin{claimproof}
			We will show that the desired upper bound holds even if we allow an adversary to specify the outcome of the algorithm for $e_0 = e$. So, fix the outcome of the algorithm's execution for $e_0 = e$ and let $\mathbb{E}^\ast$ denote conditional expectation subject to this outcome. We may assume that $e$ is lucky, since otherwise $\deg_\Gamma(e) = 0$. Then
			\begin{equation}\label{eq:expectation}
			\mathbb{E}^\ast \left[\deg_\Gamma(e)\right] \,\leq\, \sum_{u \in V(H^e)} \mathbb{E}\left[\,\left|\left\{h \in U \setminus \set{e}  \,:\, h \in L \text{ and } u \in V(H^h)\right\}\right|\,\right].
			\end{equation}
			We may write ``$\mathbb{E}$'' instead of ``$\mathbb{E}^\ast$'' on the right-hand side of \eqref{eq:expectation} because the construction of $H^h$ for $h \neq e$ is independent from the construction of $H^e$. Fix an arbitrary vertex $u \in V(H^e)$. Observe that if $u \in V(H^h)$ for some $h$, then there is some $i$ such that $u \in V(F_i^h)$ or $u \in V(P_i^h)$, so
			\[
			\mathbb{E}\left[\,\left|\left\{h \in U \setminus \set{e}  \,:\, h \in L \text{ and } u \in V(H^h)\right\}\right|\,\right] \,\leq\, \sum_{i = 0}^{\infty} \mathbb{E}\left[\,\left|\left\{h \in U \,:\, u \in V(F_i^h) \text{ or } u \in V(P_i^h)\right\}\right|\,\right].
			\]
			If $u \in V(F_i^h)$ or $u \in V(P_i^h)$, then, by Lemma~\ref{lemma:back_orbit},
			\begin{equation}\label{eq:backagain}
			x_i^h \in R^\shortleftarrow(u, \phi), \quad x_{i-1}^h \in R^\shortleftarrow(x_i^h, \phi), \quad \ldots, \quad x_1^h \in R^\shortleftarrow(x_2^h, \phi), \quad \text{and} \quad x^h = x_0^h  \in R^\shortleftarrow(x_1^h, \phi).
			\end{equation}
			Since for every $x \in V$, there are at most $\Delta$ edges $h$ such that $x^h = x$, we conclude that, by Lemma~\ref{lemma:indegree}, there are at most $\Delta(\Delta+1)^{3(i+1)} \leq (\Delta+1)^{3i + 4}$ possible sequences $(h, x_0^h, x_1^h, \ldots,x_i^h)$ satisfying \eqref{eq:backagain}. Since the vertices $x_1^h$, \ldots, $x_i^h$ are chosen uniformly at random from sets of $\ell$ candidates, we obtain
			\[
			\mathbb{E}\left[\,\left|\left\{h \in U \,:\, u \in V(F_i^h) \text{ or } u \in V(P_i^h)\right\}\right|\,\right] \,\leq\, \frac{(\Delta+1)^{3i+4}}{\ell^i} \,=\, \frac{(\Delta+1)^{4}}{\lambda^i}.
			\]
			Thus we can write
			\[
			\mathbb{E}\left[\,\left|\left\{h \in U \setminus \set{e} \,:\, h \in L \text{ and } u \in V(H^h)\right\}\right|\,\right] \,\leq\, \sum_{i=0}^\infty \frac{(\Delta+1)^{4}}{\lambda^i} \,=\, \frac{\lambda}{\lambda - 1}(\Delta+1)^4  \,=\, (1 + o(1)) (\Delta + 1)^4.
			\]
			Plugging this into \eqref{eq:expectation} and using \eqref{eq:size_bound} gives
			\[
			\mathbb{E}^\ast \left[\deg_\Gamma(e)\right] \,\leq\, (1+o(1))(\Delta+1)^4|V(H^e)| \,\leq\, (\epsilon + o(1))(\Delta+1)^{10} (\log n)^2. \qedhere
			\]
		\end{claimproof}
		
		Using the above claim and Markov's inequality, we see that for each uncolored edge $e$,
		\begin{equation}\label{eq:25}
		\P\left[\deg_\Gamma(e) > 2\epsilon(\Delta+1)^{10}(\log n)^2 - 1\right] \,\leq\, \frac{1}{2} + o(1).
		\end{equation}
		Let $W' \subseteq L$ be the set of all lucky edges $e$ satisfying $\deg_\Gamma (e) \leq 2\epsilon(\Delta+1)^{10} (\log n)^2 - 1$. Using Lemma~\ref{lemma:lucky} and \eqref{eq:25}, we conclude that for each $e \in U$,
		\begin{align*}
		\P\left[e \in W'\right] \,&\geq\, 1 \,-\, \frac{n}{\lambda^T} \,-\, \frac{3T(\Delta+1)^3}{\lambda - 1} \,-\, \P\left[\deg_\Gamma(e) > 2\epsilon(\Delta+1)^{10}(\log n)^2 - 1\right] \\
		&\geq\, 1 \,-\, n \left(\frac{1}{ (\Delta+1)^3\lfloor \log n \rfloor}\right)^{\lfloor \epsilon \log n \rfloor}\,-\, 3\epsilon \,-\, \frac{1}{2} \,-\, o(1)  \,=\, \frac{1}{2} \,-\, 3\epsilon \,-\, o(1) \,>\, \frac{1}{4}.
		\end{align*}
		Therefore, with positive probability, $|W'| \geq |U|/4$. Note that the maximum degree of the induced subgraph of $\Gamma$ on $W'$ is, by definition, at most $2\epsilon(\Delta+1)^{10} (\log n)^2 - 1$, which means that it contains an independent set $W \subseteq W'$ of size at least $|W| \geq |W'|/(2\epsilon(\Delta+1)^{10} (\log n)^2)$. Since $W$ is an independent set in $\Gamma$, the subgraphs $H^e$ for $e \in W$ are pairwise vertex-disjoint, and $|W| \geq |U|/(8\epsilon(\Delta+1)^{10}(\log n)^2) \geq |U|/((\Delta+1)^{10}(\log n)^2)$ with positive probability. This completes the proof of Theorem~\ref{theo:ext}.
	\end{scproof}
	
	\printbibliography
	
\end{document}